\documentclass{article}

\usepackage{epsfig}
\usepackage{fancybox}
\usepackage{graphicx}
\usepackage{psfrag}

\usepackage{latexsym}
\usepackage{amssymb}
\usepackage{enumerate}
\usepackage{amsmath}
\usepackage{subfigure}
\usepackage{bm}
\usepackage{hyperref}
\newcommand{\R}{{\mathbb R}}

\newcommand{\Z}{{\mathbb Z}}

\newcommand{\D}{{\partial}}

\newtheorem{thm}{Theorem}[section]

\newtheorem{cor}[thm]{Corollary}

\begin{document}

\title{Numerical resolution of an anisotropic non-linear diffusion problem
\thanks{%{Received date / Revised version date}
          % The correct dates will be entered by the CMS editor}}
}}
          %For each author, make a block with the following four macros:
\author{St\'ephane BRULL \thanks {Universit\'e de Bordeaux 1, Institut de Math\'ematiques de Bordeaux - 351, cours de la Lib\'eration, F-33405 Talence Cedex, FRANCE (email: stephane.brull@math.u-bordeaux1.fr).}
\and Fabrice DELUZET \thanks {CNRS, Institut de Math\'ematiques de Toulouse - 118 route de Narbonne, F-31062 Toulouse, FRANCE (email: fabrice.deluzet@math.univ-toulouse.fr).}
\and Alexandre MOUTON \thanks {CNRS, Laboratoire Paul Painlev\'e - Universi\'e Sciences et Technologies de Lille, Cit\'e Scientifique, F-59655 Villeneuve d'Ascq Cedex, FRANCE (email: alexandre.mouton@math.univ-lille1.fr).}}
% 
% \author{full name
% \thanks {address, (email).}
% \and full name \thanks {address, (email).}}
          %{Put the URL for your home page here if you have one}

          %Use \thanks statements for acknowledgements of grants and
          %support. They will appear below all the authors' addresses, so be
          %specific about which author is thanking whom:

          %\thanks{}

\pagestyle{myheadings} \markboth{Numerical resolution of an anisotropic non-linear diffusion problem}{S. Brull, F. Deluzet, A. Mouton}\maketitle

\begin{abstract}
This paper is devoted to the numerical resolution of an anisotropic non-linear diffusion problem involving a small parameter $\epsilon$, defined as the anisotropy strength reciprocal. In this work, the anisotropy is carried by a variable vector function $\mathbf{b}$. The equation being supplemented with Neumann boundary conditions, the limit $\epsilon \rightarrow 0$ is demonstrated to be a singular perturbation of the original diffusion equation. To address efficiently this problem, an Asymptotic-Preserving scheme is derived. This numerical method does not require the use of coordinates adapted to the anisotropy direction and exhibits an accuracy as well as a computational cost independent of the anisotropy strength.
\end{abstract}

\paragraph{keywords}
\smallskip Anisotropic diffusion problems; Singular perturbation; Asymptotic-Preserving schemes.

\paragraph{AMS Subject Classification}
\smallskip {\bf 35J60, 35J62, 65M06, 65M12, 65N06, 65N12.}

\section{Introduction}

\indent The aim of this paper is to build an efficient numerical method for solving an anisotropic diffusion problem where the anisotropy is carried by a vector $\mathbf{b}$. This work is motivated by investigations of strongly magnetized plasmas, more specifically the study of the Euler-Lorentz model in a low Mach number regime and in the presence of a large magnetic field. This framework is characteristic of the magnetically confined plasma fusion \cite{Degond,Hazeltine,Miyamoto}. In this context, the asymptotic parameter $\epsilon$ represents the gyro-period of particles as well as the square root of the Mach number, the vector field $\mathbf{b}$ being the magnetic field direction. Therefore the  $\epsilon$ values can be very small in some sub-regions of the computational domain where the magnetic field is large, inducing then a severe anisotropy of the medium, while being large in other sub-domains for intermediate and small strength of the magnetic field. Another important property of this system is the time dependence of the magnetic field defining the anisotropy direction. These two main characteristics define the framework of the present paper whose purpose is to design a numerical scheme for anisotropy ratios ranging from  $\epsilon \ll 1$ to $\epsilon\sim\mathcal{O}(1)$ and for a time varying anisotropy direction. In order to address efficiently these requirements, the numerical method should not rely on a coordinate system adapted to this anisotropy direction. The use of adapted coordinates would imply mesh modifications accordingly to the evolution of $\mathbf{b}$, an intricate and expensive procedure we wish to avoid. Thus, the numerical method introduced here will carry out the anisotropic non-linear diffusion problem on a mesh independent of the anisotropy direction.\\
\indent This scheme will be detailed on the following model problem 
\begin{equation} \label{elliptic_non-linear_eps_intro}
\left\{
\begin{array}{ll}
-\nabla_{\mathbf{x}}\cdot \left(H_{\epsilon}\,(\mathbf{b} \otimes \mathbf{b})\,\cfrac{\nabla_{\mathbf{x}} p_{\epsilon}-\mathbf{S}_{\epsilon}}{\epsilon}\right) + g_{\epsilon}(p_{\epsilon}) = f_{\epsilon} \, , & \textnormal{in $\Omega$,} \\
\left( H_{\epsilon}\,(\mathbf{b} \otimes \mathbf{b})\,\cfrac{\nabla_{\mathbf{x}} p_{\epsilon}-\mathbf{S}_{\epsilon}}{\epsilon} \right) \cdot \bm{\nu} \equiv 0 \, , & \textnormal{on $\D\Omega$.} 
\end{array}
\right.
\end{equation}
In this system, $\Omega$ is a bounded subset of $\R^{d}$ ($d = 1,2,3$), $\epsilon \geq 0$ is a fixed constant parameter and, for any $\mathbf{x} \in \D\Omega$, $\bm{\nu} = \bm{\nu}(\mathbf{x})$ stands for the unit outward normal vector. $\nabla_{\mathbf{x}}$ and $\nabla_{\mathbf{x}}\cdot$ stand for the gradient and the divergence operators with respect to the space variable $\mathbf{x}$. We assume that $H_{\epsilon} : \overline{\Omega} \to \R_{+}^{*}$, $\mathbf{b} : \overline{\Omega} \to \R^{d}-\{0\}$, $f_{\epsilon} : \Omega \to \R$, $\mathbf{S}_{\epsilon}:\overline{\Omega} \to \R^{d}$ are given and the unknown of the problem is the function $p_{\epsilon} : \overline{\Omega} \to \R$. The tensor product of two vectors $\mathbf{u}$ and $\mathbf{v}$ is denoted $\mathbf{u} \otimes \mathbf{v}$. Finally, we assume that, for any $\epsilon$, the function $p \mapsto g_{\epsilon}(p)$ is strictly increasing and can be non-linear. This equation is well suited for the plasma fusion context above depicted. It allows the computation of the plasma pressure in order to guarantee that the forces vanish in the low Mach regime for strongly magnetized plasma. The function denoted $g_\epsilon$ defines the internal energy of the fluid with respect to the pressure. This relation may be non-linear, which motivates the investigation of non-linear anisotropic problems. However, the derivation of this equation is out of the scope of the present paper and we refer to related works (see \cite{Brull-Degond-Deluzet-Mouton,Brull-Degond-Deluzet,Degond}) for detailed explanations. Furthermore, we wish to present the numerical method in a context wider than the strict plasma context, since anisotropic diffusion problem are encountered in many applications. Good examples of these applications are, for instance, image noise filtering, convection dominated diffusion equations and more generally diffusion problem with strong medium anisotropies. The model equation \eqref{elliptic_non-linear_eps_intro} is representative of a large enough variety of problems, up to slight changes, and will be considered to detail the numerical method.

\indent Developing an efficient numerical method to compute the solution of this diffusion problem, regardless to $\epsilon$ values,  is a difficult task. Indeed, the limit $\epsilon \rightarrow 0$ is a singular limit for the problem (\ref{elliptic_non-linear_eps_intro}), the   diffusion equation degenerating into the following one
\begin{equation} \label{elliptic_non-linear_0_intro}
\left\{
\begin{array}{ll}
-\nabla_{\mathbf{x}}\cdot \left(H_{0}\,(\mathbf{b} \otimes \mathbf{b})\,(\nabla_{\mathbf{x}} \tilde{p}_{0}-\mathbf{S}_{0})\right) = 0 \, , & \textnormal{in $\Omega$,} \\
\left( H_{0}\,(\mathbf{b} \otimes \mathbf{b})\,(\nabla_{\mathbf{x}} \tilde{p}_{0}-\mathbf{S}_{0}) \right) \cdot \bm{\nu} \equiv 0 \, , & \textnormal{on $\D\Omega$.} 
\end{array}
\right.
\end{equation}
The system \eqref{elliptic_non-linear_0_intro} is ill-posed, its solution being non-unique. More precisely, if $\tilde{p}_{0}$ is a solution of (\ref{elliptic_non-linear_0_intro}) and $c : \overline{\Omega} \to \R$ is a function verifying $\mathbf{b} \cdot \nabla_{\mathbf{x}}c = 0$ on $\overline{\Omega}$, then $\tilde{p}_{0}+c$ defines a new solution of (\ref{elliptic_non-linear_0_intro}). However, $p_0$ the limit of $p_\epsilon$ solution of (\ref{elliptic_non-linear_eps_intro}) is uniquely defined by the limit problem as demonstrated in Section~\ref{dec}, but a direct discretization of the diffusion problem \eqref{elliptic_non-linear_eps_intro} gives rise to a linear system with a conditioning number that blows up for vanishing $\epsilon$. This property has been outlined for numerical studies of elliptic equation singular perturbations (see \cite{Chang,Degond-Deluzet-Negulescu}). \\
\indent To tackle this difficulty, an \textit{Asymptotic-Preserving} (AP) scheme is introduced to compute the solution of the anisotropic diffusion problem for $\epsilon = \mathcal{O}(1)$ and to capture $p_0$, the solution of the limit problem, for small $\epsilon$ values. This property should be provided without any limitations on the discretization parameters related to the value of $\epsilon$. These requirements are compliant with the properties of AP-schemes originally introduced in \cite{Jin} and developed in \cite{Lemou-Mieussens} for diffusive regimes of transport equations. These techniques have received numerous extensions to other singular perturbation problems: relaxation limits of kinetic plasma descriptions \cite{Crouseilles-Lemou,Filbet-Jin_stiff,Filbet-Jin_BGK}, quasi-neutral limit of fluid and kinetic plasma models \cite{Belaouar-Crouseilles-Degond-Sonnendrucker,crisp1,crisp2,crisp3,TSAPS,Navoret,Savelief,Degond-Liu-Vignal}, hydrodynamic low Mach number limit 
 \cite{Degond-Tang,Klar}, radiative hydrodynamics \cite{Buet-Cordier-Lucquin-Mancini,Buet-Despres}, fluid and particle flows \cite{Carrillo-Goudon-Lafitte} and strongly magnetized plasmas as well as heterogeneous media \cite{Chang,Brull-Degond-Deluzet,Brull-Degond-Deluzet-Mouton,Narski,Narski2,Degond,Degond-Deluzet-Negulescu,Degond-Deluzet-Sangam-Vignal}.

\indent The Asymptotic-Preserving property of the presented method is obtained thanks to a decomposition of the solution, introduced in \cite{Degond-Deluzet-Negulescu} and also used in \cite{Chang,Brull-Degond-Deluzet,Narski}. It consists of the following identity $p_{\epsilon}=\pi_{\epsilon}+q_\epsilon$, $\pi_{\epsilon}$ being the solution mean part, with respect to the anisotropy ($\mathbf{b}$) direction, $q_\epsilon$ the fluctuating part. These two components verify $\pi_{\epsilon} \in K$ and $q_\epsilon\in K^\perp$, $K$ defining the functions constant along the  $\mathbf{b}$-direction, $K^\perp$ the functions of zero mean value along $\mathbf{b}$. This decomposition was first developed for meshes adapted to the anisotropy direction \cite{Chang,Degond-Deluzet-Negulescu}, for which, the discretization of $K$ is straightforward. A direct discretization of the sub-space $K^\perp$ is, on the other side, much more intricate. This difficulty is overcome thanks to the introduction of a Lagrangian multiplier, in order to penalize the zero mean value property of the functions belonging to $K^\perp$. The method is extended in \cite{Narski} for computations with meshes independent of the anisotropy direction. This is achieved by introducing two more Lagrangian multipliers to discretize the sub-spaces. The size of the linear system providing the problem solution is then significantly enlarged. However, this drawback may be corrected thanks to a slightly different decomposition. In \cite{Narski2} the solution is decomposed in two non-orthogonal parts which allows  the definition of two sub-spaces whose direct discretization is readily obtained without any Lagrangian multipliers. The size of the linear system obtained with this approach is considerably lowered compared to the previous method \cite{Narski}. This method has been extended in \cite{MN2012,LNN2013} to non-linear diffusion equations. The path followed in the present paper still relies on the decomposition in $K$ and $K^\perp$. However, the discretization of these sub-spaces is achieved using a differential characterization, similar to the one introduced in \cite{Brull-Degond-Deluzet}. This finally allows the computation of the solution thanks to a second-order problem for $\pi_{\epsilon}$ and a fourth-order problem for $q_{\epsilon}$. This latter problem, in the framework of Neumann boundary conditions considered in this paper, can be recast into two elliptic problems.

\indent The method proposed finally reduces to the computation of three standard elliptic problems for which very efficient solvers can be used (for instance multi-grid solvers). For the former approaches \cite{Narski,Narski2}, the equations providing both components are not classical elliptic equations and the resolution of the linear system requires more sophisticated solvers. This complexity is resource demanding and may be challenging for realistic three-dimensional computations. Finally this paper also presents an extension to non-linear reaction diffusion problems, a class of problems that has never been investigated in the previous works.

\indent The paper is organized as follows: in Section 2, the decomposition methodology is presented. The linear case, \textit{i.e.} with $g_{\epsilon}(p)(\mathbf{x}) = G_{\epsilon}(\mathbf{x})\,p(\mathbf{x})$ where $G_{\epsilon} : \overline{\Omega} \to \R_{+}^{*}$ is a given function sequence, is first investigated: more precisely, we describe the decomposition procedure in the specific case where $G_{\epsilon}$ is a strictly positive constant denoted $\lambda_{\epsilon}$, then we generalize this procedure to any function $G_{\epsilon} : \overline{\Omega} \to \R_{+}^{*}$ by using well-chosen Sobolev spaces. Finally the non-linear problems are addressed by invoking Gummel's iterative algorithm. Section 3 is devoted to presentation of the discretization. Finally, the efficiency of the numerical method is demonstrated in Section 4.

\section{Scale separation and solution decomposition}
\label{dec}
\setcounter{equation}{0}

\indent In this section a scale separation is introduced to ensure the Asymptotic-Preserving property of the scheme. This is achieved by transforming the singular perturbation problem \eqref{elliptic_non-linear_eps_intro} into an equivalent system for which the limit $\epsilon \rightarrow 0$ is regular. For simplicity reasons, the linear case with constant $G_{\epsilon}$ is first considered for detailing the decomposition method. In this framework, the singular nature of the limit $\epsilon\rightarrow 0$ is outlined and the limit problem, providing $p_0=\lim_{\epsilon\rightarrow 0} p_\epsilon$, is stated. A development to linear cases with variable positive functions $G_{\epsilon}$ is then presented and finally, thanks to Gummel's iterative method \cite{Gummel}, the non-linear case is addressed by using a sequence of linear problems.

\subsection{AP-scheme derivation for linear problems}

\subsubsection{A simplified framework: constant $G_{\epsilon}$}

We assume here that the given sequence $(g_{\epsilon})_{\epsilon \, \geq \, 0}$ is of the form
\begin{equation*}
g_{\epsilon}(p)(\mathbf{x}) = \lambda_{\epsilon} \, p(\mathbf{x}) \, ,
\end{equation*}
where $\lambda_{\epsilon} > 0$ is a known constant for any $\epsilon \geq 0$. Then the diffusion problem (\ref{elliptic_non-linear_eps_intro}) writes
\begin{equation} \label{elliptic_linear_eps}
\left\{
\begin{array}{ll}
-\nabla_{\mathbf{x}}\cdot \left(H_{\epsilon}\,(\mathbf{b} \otimes \mathbf{b})\,(\nabla_{\mathbf{x}} p_{\epsilon}-\mathbf{S}_{\epsilon})\right) + \epsilon\, \lambda_{\epsilon} \, p_{\epsilon} = \epsilon \, f_{\epsilon} \, , & \textnormal{in $\Omega$,} \\
\left( H_{\epsilon}\,(\mathbf{b} \otimes \mathbf{b})\,(\nabla_{\mathbf{x}} p_{\epsilon}-\mathbf{S}_{\epsilon}) \right) \cdot \bm{\nu} \equiv 0 \, , & \textnormal{on $\D\Omega$\,.} 
\end{array}
\right.
\end{equation}
 The limit solution $p_{0}$ of the singular perturbation problem \eqref{elliptic_linear_eps} verifies the limit problem 
\begin{equation}\label{eq:def:limit:problem}
\left\{
\begin{array}{ll}
\displaystyle -\lim_{\epsilon\,\to\,0} \nabla_{\mathbf{x}} \cdot \left( H_{\epsilon}\,(\mathbf{b} \otimes \mathbf{b})\,\cfrac{\nabla_{\mathbf{x}}p_{\epsilon}-\mathbf{S}_{\epsilon}}{\epsilon}\right) + \lambda_{0}\,p_{0} = f_{0} \, , & \textnormal{in $\Omega$,} \\
\displaystyle \left( H_{0}\,(\mathbf{b} \otimes \mathbf{b})\,(\nabla_{\mathbf{x}} p_{0}-\mathbf{S}_{0}) \right) \cdot \bm{\nu} \equiv 0 \, , & \textnormal{on $\D\Omega$.}
\end{array}
\right.
\end{equation}
This algebraic equation admits a unique solution under the assumption 
\begin{equation}\label{eq:limit:requirement}
H_{\epsilon}\,(\mathbf{b} \otimes \mathbf{b})(\nabla_{\mathbf{x}}p_{\epsilon}-\mathbf{S}_{\epsilon}) = \mathcal{O}(\epsilon) \, ,
\end{equation}
a requirement that must be fulfilled by the numerical method. To ensure this property, the methodology consists in using a decomposition similar to that of \cite{Brull-Degond-Deluzet,Brull-Degond-Deluzet-Mouton,TSAPS,Narski,Degond-Deluzet-Negulescu}. The solution is decomposed into $\pi_\epsilon$ its mean part with respect to the anisotropy direction and $q_\epsilon$ the fluctuating part, which exhibits the property to have a zero mean value along the anisotropy direction. These two functions verify $\pi_\epsilon \in K$ and $q_\epsilon \in K^\perp$, $K$ being the kernel of the elliptic operator defined by equation~\ref{elliptic_non-linear_0_intro}. These properties are capitalized on, to isolate in the problem \ref{elliptic_linear_eps} the macro scale (providing $\pi_\epsilon$) from the micro scale (giving $q_\epsilon$) and thereby, build the Asymptotic-Preserving scheme. The main difficulty of the procedure lies in the characterization of the sub-spaces associated to the different scales. In \cite{Chang,TSAPS,Narski,Degond-Deluzet-Negulescu} the property of the functions populating $K$ or $K^\perp$ are imposed by a penalization technique. The methodology developed in this paper operates a similar decomposition on to $K$ and $K^\perp$, but with a different characterization of these sub-spaces. Here, we shape the technique introduced in \cite{Brull-Degond-Deluzet} for a very specific framework, in order to discriminate the functions in $K$ and $K^\perp$ thanks to differential properties, providing thus, an easy discretization. 

With this aim, we introduce the following Sobolev spaces:
\begin{equation*}
\begin{split}
V &= \left\{ p \in L^{2}(\Omega) \, : \, \mathbf{b}\cdot \nabla_{\mathbf{x}}p \in L^{2}(\Omega) \, \right\} \, , \\
W &= \left\{ q \in L^{2}(\Omega) \, : \, \nabla_{\mathbf{x}} \cdot (\mathbf{b} \, q) \in L^{2}(\Omega) \, \right\} \, , \\
W_{0} &= \left\{ q \in W \, : \, (\mathbf{b}\,q) \cdot \bm{\nu} \equiv 0 \, \, \textnormal{on $\D\Omega$} \, \right\} \, , \\
\end{split}
\end{equation*}
and we define $K \subset V$ as
\begin{equation*}
K = \left\{ \pi \in V \, : \, \mathbf{b}\cdot \nabla_{\mathbf{x}} \pi = 0 \, \, \textnormal{on $\Omega$} \, \right\} \, .
\end{equation*}
The goal is to reproduce the function decomposition into its mean and fluctuating parts. The functions of $K$ correspond to the mean part and the complementary part is demonstrated to belong to $K^\perp$. This is the purpose of the following theorem:
\begin{thm} \label{theorem_decompo_linear}
We denote by $\nabla_{\mathbf{x}} \cdot (\mathbf{b} \, W_{0})$ the subspace of functions $\theta \in L^{2}(\Omega)$ such that
\begin{equation}
\exists \, \chi \in W_{0} \, , \qquad \theta = \nabla_{\mathbf{x}} \cdot (\mathbf{b}\,\chi) \, ,
\end{equation}
and we equip it with the usual norm on $L^{2}(\Omega)$. Then:
\begin{itemize}
\item $W_{0}$ equipped with the norm $\|p\|_{W_{0}} = \left\|\nabla_{\mathbf{x}} \cdot (\mathbf{b}\,p) \right\|_{L^{2}}$ is a Hilbert space,
\item $\nabla_{\mathbf{x}} \cdot (\mathbf{b} \, W_{0})$ is a closed subspace in $L^{2}(\Omega)$,
\item $K$ is a closed subspace in $L^{2}(\Omega)$.
\item We have the orthogonal decomposition
\begin{equation} \label{orthogonal_decomposition_linear}
L^{2}(\Omega) = K \oplus K^{\perp} \, , \quad \textit{with} \quad K^{\perp} = \nabla_{\mathbf{x}} \cdot (\mathbf{b} \, W_{0}) \, .
\end{equation}
\end{itemize}
\end{thm}
The demonstration of this theorem will be omitted. It can be readily adapted from that of Theorem 2.1 from \cite{Brull-Degond-Deluzet}. As a consequence of this theorem, the decomposition 
\begin{equation} \label{decompo_linear}
p_{\epsilon} = \pi_{\epsilon} + q_{\epsilon} \, , \, \pi_{\epsilon} \in K, \, q_{\epsilon} \in K^{\perp} \, ,
\end{equation}
exists and is unique for any $\epsilon \geq 0$. Therefore finding the particular solution $p_{0}$ which is exactly the limit of $(p_{\epsilon})_{\epsilon\,>\,0}$ is equivalent to find $\pi_{0}$ and $q_{0}$ as the respective limits of $(\pi_{\epsilon})_{\epsilon\,>\,0}$ and $(q_{\epsilon})_{\epsilon\,>\,0}$. Then, our goal is now to find some equations for $\pi_{\epsilon}$ and $q_{\epsilon}$ which are well-posed for any value of $\epsilon$, including $\epsilon = 0$. For this purpose, the decomposition (\ref{decompo_linear}) is introduced into (\ref{elliptic_linear_eps}), yielding
\begin{equation} \label{elliptic_linear_eps_decomposed}
\left\{
\begin{array}{ll}
-\nabla_{\mathbf{x}}\cdot \left(H_{\epsilon}\,(\mathbf{b} \otimes \mathbf{b})\,(\nabla_{\mathbf{x}} q_{\epsilon}-\mathbf{S}_{\epsilon})\right) + \epsilon\, \lambda_{\epsilon} \, (\pi_{\epsilon}+q_{\epsilon}) = \epsilon \, f_{\epsilon} \, , & \textnormal{in $\Omega$,} \\
\left( H_{\epsilon}\,(\mathbf{b} \otimes \mathbf{b})\,(\nabla_{\mathbf{x}} q_{\epsilon}-\mathbf{S}_{\epsilon}) \right) \cdot \bm{\nu} \equiv 0 \, , & \textnormal{on $\D\Omega$.} 
\end{array}
\right.
\end{equation}
The variational formulation on $V$ writes
\begin{equation} \label{weak_linear_eps}
\begin{split}
\int_{\Omega} H_{\epsilon} \, (\mathbf{b} \cdot \nabla_{\mathbf{x}} q_{\epsilon}) \, (\mathbf{b} \cdot \nabla_{\mathbf{x}} \theta) \, d\mathbf{x} + \epsilon \, \lambda_{\epsilon} \int_{\Omega} (\pi_{\epsilon} + q_{\epsilon}) \, \theta \, d\mathbf{x} \\
= \epsilon \int_{\Omega} f_{\epsilon} \, \theta \, dx + \int_{\Omega} H_{\epsilon}\, (\mathbf{b} \cdot \mathbf{S}_{\epsilon}) (\mathbf{b} \cdot \nabla_{\mathbf{x}}\theta) \, d\mathbf{x} \, ,
\end{split}
\end{equation}
for any test function $\theta \in V$. \\

\indent In order to exhibit the equation providing $\pi_{\epsilon} \in K$, the variational formulation (\ref{weak_linear_eps}) is tested against $\theta \in K$ giving
\begin{equation*}
\int_{\Omega} (\lambda_{\epsilon} \, \pi_{\epsilon} - f_{\epsilon}) \, \theta \, d\mathbf{x} = 0 \, ,
\end{equation*}
which means that $\lambda_{\epsilon}\,\pi_{\epsilon}-f_{\epsilon} \in K^{\perp}$ for any $\epsilon \geq 0$. According to Theorem \ref{theorem_decompo_linear}, there exists a function $h_{\epsilon} \in W_{0}$ such that
\begin{equation} \label{eq}
\lambda_{\epsilon} \, \pi_{\epsilon} - f_{\epsilon} = \nabla_{\mathbf{x}} \cdot (\mathbf{b} \, h_{\epsilon}) \, .
\end{equation}
This equation furnishes a means of computation for $\pi_\epsilon$. Firstly, applying the differential operator $\mathbf{b} \cdot \nabla_{\mathbf{x}}$ onto (\ref{eq}) leads to an equation for $h_\epsilon$
\begin{equation} \label{problem_heps_linear}
\left\{
\begin{array}{ll}
-\mathbf{b} \cdot \nabla_{\mathbf{x}} \left(\nabla_{\mathbf{x}} \cdot (\mathbf{b} \, h_{\epsilon}) \right) = \mathbf{b} \cdot \nabla_{\mathbf{x}} f_{\epsilon} \, , & \textnormal{in $\Omega$,} \\
(\mathbf{b} \, h_{\epsilon}) \cdot \bm{\nu} \equiv 0 \, , & \textnormal{on $\D\Omega$,}
\end{array}
\right.
\end{equation}
then, $\pi_\epsilon$ is retrieved thanks to
\begin{equation} \label{def_pieps_linear}
\pi_{\epsilon} = \cfrac{1}{\lambda_{\epsilon}} \, \left[ f_{\epsilon} + \nabla_{\mathbf{x}} \cdot (\mathbf{b} \, h_{\epsilon})\right] \, .
\end{equation}
Note that the system (\ref{problem_heps_linear})-(\ref{def_pieps_linear}) is well-posed  and does not degenerate for any value of $\epsilon \geq 0$, including $\epsilon = 0$. It provides a means of computing the macro component of the solution regardless to $\epsilon$ values. \\

\indent To derive an equation for $q_{\epsilon} \in K^{\perp}$, we now assume that the test function $\theta$ in (\ref{weak_linear_eps}) is in $K^{\perp}$. According to Theorem \ref{theorem_decompo_linear}, there exist two functions $\chi$ and $l_{\epsilon}$ in $W_{0}$ such that
\begin{equation} \label{def_qeps_linear}
q_{\epsilon} = \nabla_{\mathbf{x}} \cdot (\mathbf{b} \, l_{\epsilon}) \, ,
\end{equation}
and
\begin{equation*}
\theta = \nabla_{\mathbf{x}} \cdot (\mathbf{b} \, \chi) \, .
\end{equation*}
As a consequence, the variational formulation of (\ref{elliptic_linear_eps_decomposed}) can be rewritten as follows:
\begin{equation*}
\begin{split}
\int_{\Omega} H_{\epsilon} \left(\mathbf{b} \cdot \nabla_{\mathbf{x}} \left( \nabla_{\mathbf{x}} \cdot (\mathbf{b} \, l_{\epsilon}) \right) \right) \, \left(\mathbf{b} \cdot \nabla_{\mathbf{x}} \left( \nabla_{\mathbf{x}} \cdot (\mathbf{b} \, \chi) \right) \right) \, d\mathbf{x} + \epsilon\, \lambda_{\epsilon} \int_{\Omega}  \left( \nabla_{\mathbf{x}} \cdot (\mathbf{b} \, l_{\epsilon}) \right) \, \left(\nabla_{\mathbf{x}} \cdot (\mathbf{b} \, \chi) \right) \, d\mathbf{x} \\
= \epsilon \int_{\Omega} f_{\epsilon} \, \left(\nabla_{\mathbf{x}} \cdot (\mathbf{b} \, \chi) \right) \, d\mathbf{x} + \int_{\Omega} H_{\epsilon}\, (\mathbf{b} \cdot \mathbf{S}_{\epsilon}) \left(\mathbf{b} \cdot \nabla_{\mathbf{x}} \left( \nabla_{\mathbf{x}} \cdot (\mathbf{b}\,\chi)\right) \right) \, d\mathbf{x} \, .
\end{split}
\end{equation*}
We recognize the variational formulation of
\begin{equation} \label{problem_leps_linear}
\hspace{-0.2cm}\left\{
\begin{array}{ll}
\mathbf{b} \cdot \nabla_{\mathbf{x}} \left( \nabla_{\mathbf{x}} \cdot \left(H_{\epsilon} \, (\mathbf{b} \otimes \mathbf{b}) \, \nabla_{\mathbf{x}} \left( \nabla_{\mathbf{x}} \cdot (\mathbf{b} \, l_{\epsilon}) \right) \right) \right) - \epsilon \, \lambda_{\epsilon}\, \mathbf{b} \cdot \nabla_{\mathbf{x}} \left( \nabla_{\mathbf{x}} \cdot (\mathbf{b} \, l_{\epsilon}) \right) & \\
\begin{array}{r}
\qquad \qquad \qquad \qquad \qquad = - \mathbf{b} \cdot \nabla_{\mathbf{x}} \left( \epsilon\, f_{\epsilon} - \nabla_{\mathbf{x}} \cdot \left(H_{\epsilon}\,(\mathbf{b} \otimes \mathbf{b})\,\mathbf{S}_{\epsilon} \right) \right) \, ,
\end{array}
& \textnormal{in $\Omega$,} \\
\left(H_{\epsilon} \, (\mathbf{b} \otimes \mathbf{b}) \, \nabla_{\mathbf{x}} \left( \nabla_{\mathbf{x}} \cdot (\mathbf{b} \, l_{\epsilon}) \right) \right) \cdot \bm{\nu} \equiv \left(H_{\epsilon}\, (\mathbf{b} \otimes \mathbf{b}) \, \mathbf{S}_{\epsilon}\right) \cdot \bm{\nu} \, , & \textnormal{on $\D\Omega$,} \\
(\mathbf{b} \, l_{\epsilon}) \cdot \bm{\nu} \equiv 0 \, , & \textnormal{on $\D\Omega$.}
\end{array}
\right.
\end{equation}
Therefore, coupling this system with (\ref{def_qeps_linear}), we recognize a complete definition of $q_{\epsilon}$ which is well-posed for any $\epsilon \geq 0$, including $\epsilon = 0$. Moreover this computation of $q_\epsilon$ is totally compliant with the condition~\eqref{eq:limit:requirement} and guarantees the Asymptotic-Preserving property of the scheme.\\

\indent At this point, we have established a system of equations for $\pi_{\epsilon}$ and $q_{\epsilon}$ which is well-posed for any $\epsilon > 0$ but also for $\epsilon = 0$. Then, solving the well-posed equations (\ref{problem_heps_linear}), (\ref{problem_leps_linear}), (\ref{def_pieps_linear}) and (\ref{def_qeps_linear}) provides $\pi_{0}$ and $q_{0}$ as the respective limits of $(\pi_{\epsilon})_{\epsilon\,>\,0}$ and $(q_{\epsilon})_{\epsilon\,>\,0}$ when $\epsilon \to 0$. As a consequence, the sum $\pi_{0}+q_{0}$ is exactly the solution $p_{0}$ of (\ref{eq:def:limit:problem}). Furthermore, we can remark that the limit $\epsilon \to 0$ is regular for the reformulated model (\ref{problem_heps_linear})-(\ref{problem_leps_linear})-(\ref{def_pieps_linear})-(\ref{def_qeps_linear}).

\subsubsection{Case with variable $G_{\epsilon}$}

\indent In this paragraph, we extend the method we have presented to the general linear case, \textit{i.e.} to cases where $g_{\epsilon}$ is of the form
\begin{equation*}
g_{\epsilon}(p)(\mathbf{x}) = G_{\epsilon}(\mathbf{x}) \, p(\mathbf{x}) \, .
\end{equation*}
$G_{\epsilon} : \overline{\Omega} \to \R$ is given for any $\epsilon \geq 0$, and is supposed to be strictly positive on $\overline{\Omega}$. In such a case, the diffusion problem (\ref{elliptic_non-linear_eps_intro}) writes
\begin{equation} \label{elliptic_quasi-linear_eps}
\left\{
\begin{array}{ll}
-\nabla_{\mathbf{x}}\cdot \left(H_{\epsilon}\,(\mathbf{b} \otimes \mathbf{b})\,(\nabla_{\mathbf{x}} p_{\epsilon}-\mathbf{S}_{\epsilon})\right) + \epsilon\, G_{\epsilon} \, p_{\epsilon} = \epsilon \, f_{\epsilon} \, , & \textnormal{in $\Omega$,} \\
\left( H_{\epsilon}\,(\mathbf{b} \otimes \mathbf{b})\,(\nabla_{\mathbf{x}} p_{\epsilon}-\mathbf{S}_{\epsilon}) \right) \cdot \bm{\nu} \equiv 0 \, , & \textnormal{on $\D\Omega$.} 
\end{array}
\right.
\end{equation}
The study of these cases is motivated by the fact that the use of Gummel's algorithm on the non-linear case leads to the resolution of a sequence of linearized problems which are similar to (\ref{elliptic_quasi-linear_eps}). We refer to Section \ref{Gummel_algo} for more details about the linearization procedure. \\

\indent In order to solve the linear problem (\ref{elliptic_quasi-linear_eps}) for any value of $\epsilon$, we use the method presented in the previous paragraph. Firstly, we define $L^{2}(\Omega; G_{\epsilon})$ by
\begin{equation*}
L^{2}(\Omega; G_{\epsilon}) = \left\{ p : \Omega \to \R, \; \|p\|_{L^{2}(\Omega; G_{\epsilon})}^{2} = \int_{\Omega} G_{\epsilon}(\mathbf{x}) \, \left|p(\mathbf{x})\right|^{2} \, d\mathbf{x} < +\infty \right\}\, .
\end{equation*}
Then we introduce the following weighted Sobolev spaces:
\begin{equation*}
\begin{split}
V_{\epsilon} &= \left\{ p \in L^{2}(\Omega; G_{\epsilon}) \, : \, \mathbf{b}\cdot \nabla_{\mathbf{x}}p \in L^{2}(\Omega; G_{\epsilon}) \, \right\} \, , \\
W_{\epsilon} &= \left\{ q \in L^{2}(\Omega; G_{\epsilon}) \, : \, \nabla_{\mathbf{x}} \cdot (G_{\epsilon}\,\mathbf{b} \, q) \in L^{2}(\Omega; G_{\epsilon}) \, \right\} \, , \\
W_{0,\epsilon} &= \left\{ q \in W_{\epsilon} \, : \,
(G_{\epsilon}\,\mathbf{b}\,q) \cdot \bm{\nu} \equiv 0 \, \,
\textnormal{on $\D\Omega$} \, \right\} \, , 
\end{split}
\end{equation*}
and the set representing the functions constant along the magnetic
field lines
\begin{equation*}
K_{\epsilon} = \left\{ \pi \in V_{\epsilon} \, : \, \mathbf{b}\cdot \nabla_{\mathbf{x}} \pi = 0 \, \, \textnormal{on $\Omega$} \, \right\} \, .
\end{equation*}
Following the methodology presented in the previous paragraph and in \cite{Brull-Degond-Deluzet}, we deduce
\begin{cor} \label{theorem_decompo_quasi-linear}
$W_{0,\epsilon}$ equipped with the norm $\|p\|_{W_{0,\epsilon}} = \left\|\nabla_{\mathbf{x}} \cdot (G_{\epsilon}\,\mathbf{b}\,p) \right\|_{L^{2}(\Omega;G_{\epsilon})}$ is a Hilbert space and $\nabla_{\mathbf{x}} \cdot (G_{\epsilon}\,\mathbf{b} \, W_{0,\epsilon})$ is a closed space in $L^{2}(\Omega; G_{\epsilon})$. Furthermore, $K_{\epsilon}$ is also a closed space in $L^{2}(\Omega; G_{\epsilon})$ and we have the orthogonal decomposition
\begin{equation} \label{orthogonal_decomposition_quasi-linear}
L^{2}(\Omega; G_{\epsilon}) = K_{\epsilon} \oplus K_{\epsilon}^{\perp} \, , \quad \textit{with} \quad K_{\epsilon}^{\perp} = \cfrac{1}{G_{\epsilon}}\, \nabla_{\mathbf{x}} \cdot (G_{\epsilon}\, \mathbf{b} \, W_{0,\epsilon}) \, .
\end{equation}
\end{cor}

From the orthogonal decomposition (\ref{orthogonal_decomposition_quasi-linear}), the solution $p_{\epsilon}$ of (\ref{elliptic_quasi-linear_eps}) can be uniquely decomposed as
\begin{equation} \label{decompo_p_quasi-linear}
p_{\epsilon} = \pi_{\epsilon} + q_{\epsilon} \, , \, \pi_{\epsilon} \in K_{\epsilon}, \, q_{\epsilon} \in K_{\epsilon}^{\perp}\, .
\end{equation}
Then, if we identify the limits $\pi_{0}$ and $q_{0}$ of the sequences $(\pi_{\epsilon})_{\epsilon\,>\,0}$ and $(q_{\epsilon})_{\epsilon\,>\,0}$, we will find the limit $p_{0}$ of $(p_{\epsilon})_{\epsilon\,>\,0}$ by taking $p_{0} = \pi_{0} + q_{0}$.

In order to identify a set of equations satisfied by $\pi_{\epsilon}$ and $q_{\epsilon}$, we follow the same procedure as in the previous paragraph: we multiply (\ref{elliptic_quasi-linear_eps}) by a test function $\theta \in V_{\epsilon}$ and we integrate over $\Omega$. By choosing $\theta$ in $K_{\epsilon}$ or in $K_{\epsilon}^{\perp}$, we prove that $\pi_{\epsilon}$ and $q_{\epsilon}$ are respectively of the form
\begin{equation} \label{quasi-linear_pieps_qeps}
\begin{split}
\pi_{\epsilon} = \cfrac{1}{G_{\epsilon}} \, \left[ f_{\epsilon} + \nabla_{\mathbf{x}} \cdot (G_{\epsilon} \, \mathbf{b} \, h_{\epsilon}) \right] \, , \;
q_{\epsilon} = \cfrac{1}{G_{\epsilon}} \, \nabla_{\mathbf{x}} \cdot (G_{\epsilon} \, \mathbf{b} \, l_{\epsilon}) \, , 
\end{split}
\end{equation}
where $h_{\epsilon}$ and $l_{\epsilon}$ are solutions of
\begin{equation} \label{quasi-linear_h}
\left\{
\begin{array}{ll}
-\mathbf{b} \cdot \nabla_{\mathbf{x}} \left(\cfrac{1}{G_{\epsilon}} \, \nabla_{\mathbf{x}} \cdot (G_{\epsilon} \, \mathbf{b} \, h_{\epsilon}) \right) = \mathbf{b} \cdot \nabla_{\mathbf{x}} \left( \cfrac{f_{\epsilon}}{G_{\epsilon}} \right) \, , & \textnormal{in $\Omega$,} \\
(G_{\epsilon} \, \mathbf{b} \, h_{\epsilon}) \cdot \bm{\nu} \equiv 0 \, , & \textnormal{on $\D\Omega$,}
\end{array}
\right.
\end{equation}
and
\begin{equation} \label{quasi-linear_l}
\left\{
\begin{array}{ll}
\mathbf{b} \cdot \nabla_{\mathbf{x}} \left( \cfrac{1}{G_{\epsilon}} \, \nabla_{\mathbf{x}} \cdot \left(H_{\epsilon} \, (\mathbf{b} \otimes \mathbf{b}) \, \nabla_{\mathbf{x}} \left( \cfrac{1}{G_{\epsilon}} \, \nabla_{\mathbf{x}} \cdot (G_{\epsilon} \, \mathbf{b} \, l_{\epsilon}) \right) \right) \right) & \\
- \epsilon \, \mathbf{b} \cdot \nabla_{\mathbf{x}} \left( \cfrac{1}{G_{\epsilon}} \, \nabla_{\mathbf{x}} \cdot (G_{\epsilon} \, \mathbf{b} \, l_{\epsilon}) \right) \\
\qquad = - \mathbf{b} \cdot \nabla_{\mathbf{x}} \left( \cfrac{1}{G_{\epsilon}}\, \left(\epsilon\, f_{\epsilon} - \nabla_{\mathbf{x}} \cdot \left(H_{\epsilon}\,(\mathbf{b}\otimes \mathbf{b})\, \mathbf{S}_{\epsilon}\right) \right) \right) \, , & \textnormal{in $\Omega$,} \\
\left[H_{\epsilon} \, (\mathbf{b} \otimes \mathbf{b}) \, \nabla_{\mathbf{x}} \left( \cfrac{1}{G_{\epsilon}} \, \nabla_{\mathbf{x}} \cdot (G_{\epsilon} \, \mathbf{b} \, l_{\epsilon}) \right) \right] \cdot \bm{\nu} \equiv \left(H_{\epsilon}\,(\mathbf{b}\otimes\mathbf{b})\,\mathbf{S}_{\epsilon}\right) \cdot \bm{\nu} \, , & \textnormal{on $\D\Omega$,} \\
(G_{\epsilon} \, \mathbf{b} \, l_{\epsilon}) \cdot \bm{\nu} \equiv 0 \, , & \textnormal{on $\D\Omega$.}
\end{array}
\right.
\end{equation}
As in the previous paragraph, we observe that the equations (\ref{quasi-linear_pieps_qeps})-(\ref{quasi-linear_h})-(\ref{quasi-linear_l}) remain well-posed for any $\epsilon \geq 0$. As a consequence, the particular solution $p_{0}$ of the limit problem we are looking for is exactly the sum $\pi_{0}+q_{0}$ where $\pi_{0}$ and $q_{0}$ are computed by solving (\ref{quasi-linear_pieps_qeps})-(\ref{quasi-linear_h})-(\ref{quasi-linear_l}) with $\epsilon = 0$. \\
\indent Furthermore, the resolution of the fourth order problem (\ref{quasi-linear_l}) can be replaced by the successive resolution of two homogeneous Dirichlet type problems which are
\begin{equation} \label{eq_L_quasi-linear_eps}
\left\{
\begin{array}{ll}
-\mathbf{b} \cdot \nabla_{\mathbf{x}} \left( \cfrac{1}{G_{\epsilon}} \, \nabla_{\mathbf{x}} \cdot (H_{\epsilon} \, \mathbf{b} \, L_{\epsilon}) \right) + \epsilon \, L_{\epsilon} = - \epsilon\, \mathbf{b} \cdot \left[ \nabla_{\mathbf{x}} \left( \cfrac{f_{\epsilon}}{G_{\epsilon}}\right) -\mathbf{S}_{\epsilon} \right] \, , & \textnormal{in $\Omega$,} \\
(H_{\epsilon} \, \mathbf{b} \, L_{\epsilon}) \cdot \bm{\nu} \equiv 0 \, , & \textnormal{on $\D\Omega$,}
\end{array}
\right.
\end{equation}
and
\begin{equation} \label{eq_l_quasi-linear_eps}
\left\{
\begin{array}{ll}
-\mathbf{b} \cdot \nabla_{\mathbf{x}} \left( \cfrac{1}{G_{\epsilon}} \, \nabla_{\mathbf{x}} \cdot (G_{\epsilon} \, \mathbf{b} \, l_{\epsilon}) \right) = L_{\epsilon} -\mathbf{b} \cdot \mathbf{S}_{\epsilon} \, , & \textnormal{in $\Omega$,} \\
(G_{\epsilon} \, \mathbf{b} \, l_{\epsilon}) \cdot \bm{\nu} \equiv 0 \, , & \textnormal{on $\D\Omega$.}
\end{array}
\right.
\end{equation}

\subsection{AP-scheme derivation for non-linear problems} \label{Gummel_algo}

Finally, we consider the general model (\ref{elliptic_non-linear_eps_intro}) given in the introduction when the function $p \mapsto g_{\epsilon}(p)$ is non-linear. When $\epsilon$ goes to 0, the model becomes
\begin{equation} \label{elliptic_non-linear_0}
\left\{
\begin{array}{ll}
-\nabla_{\mathbf{x}}\cdot \left(H_{0}\,(\mathbf{b} \otimes \mathbf{b})\,(\nabla_{\mathbf{x}} \tilde{p}_{0}-\mathbf{S}_{0})\right) = 0 \, , & \textnormal{in $\Omega$,} \\
\left( H_{0}\,(\mathbf{b} \otimes \mathbf{b})\,(\nabla_{\mathbf{x}} \tilde{p}_{0}-\mathbf{S}_{0}) \right) \cdot \bm{\nu} \equiv 0 \, , & \textnormal{on $\D\Omega$.} 
\end{array}
\right.
\end{equation}
Due to the non-linearity of the function $p \mapsto g_{\epsilon}(p)$ the orthogonal decomposition method cannot be used. Then we choose to linearize the diffusion equation (\ref{elliptic_non-linear_eps_intro}) by using Gummel's algorithm developed in \cite{Gummel}. This iterative method consists in the approximation of the solution $p_{\epsilon}$ by a sequence $(p_{\epsilon,N})_{N\,\geq\,0}$ defined by
\begin{equation} \label{def_Gummel}
p_{\epsilon,N+1} = p_{\epsilon,N} + \delta_{\epsilon,N} \, ,
\end{equation}
and initialized with an arbitrary $p_{\epsilon,0}$. In this method, each $\delta_{\epsilon,N}$ is viewed as a small correction of $p_{\epsilon,N}$ in order to obtain $p_{\epsilon,N+1}$. Then, assuming that $p_{\epsilon,N+1}$ is a solution of (\ref{elliptic_non-linear_eps_intro}), it holds that
\begin{equation} \label{elliptic_non-linear_non-linearized_eps}
\left\{
\begin{array}{ll}
-\nabla_{\mathbf{x}}\cdot \left(H_{\epsilon}\,(\mathbf{b} \otimes \mathbf{b})\,\cfrac{\nabla_{\mathbf{x}} p_{\epsilon,N} + \nabla_{\mathbf{x}} \delta_{\epsilon,N} -\mathbf{S}_{\epsilon}}{\epsilon}\right) & \\
\qquad \qquad \qquad + g_{\epsilon}(p_{\epsilon,N}) + \delta_{\epsilon,N}\,g_{\epsilon}'(p_{\epsilon,N}) + \mathcal{O}(\delta_{\epsilon,N}^{2}) = f_{\epsilon} \, , & \textnormal{in $\Omega$,} \\
\left( H_{\epsilon}\,(\mathbf{b} \otimes \mathbf{b})\,\cfrac{\nabla_{\mathbf{x}} p_{\epsilon,N} + \nabla_{\mathbf{x}} \delta_{\epsilon,N}-\mathbf{S}_{\epsilon}}{\epsilon} \right) \cdot \bm{\nu} \equiv 0 \, , & \textnormal{on $\D\Omega$.} 
\end{array}
\right.
\end{equation}
Then, neglecting second order terms in $\delta_{\epsilon,N}$, we obtain a linear diffusion problem for $\delta_{\epsilon,N}$ which writes
\begin{equation} \label{elliptic_non-linear_linearized_eps}
\left\{
\begin{array}{ll}
-\nabla_{\mathbf{x}}\cdot \left(H_{\epsilon}\,(\mathbf{b} \otimes \mathbf{b})\,(\nabla_{\mathbf{x}} \delta_{\epsilon,N}-\mathbf{S}_{\epsilon,N})\right) + \epsilon\, G_{\epsilon,N} \, \delta_{\epsilon,N} = \epsilon \, f_{\epsilon,N} \, , & \textnormal{in $\Omega$,} \\
\left( H_{\epsilon}\,(\mathbf{b} \otimes \mathbf{b})\,(\nabla_{\mathbf{x}} \delta_{\epsilon,N}-\mathbf{S}_{\epsilon,N}) \right) \cdot \bm{\nu} \equiv 0 \, , & \textnormal{on $\D\Omega$,} 
\end{array}
\right.
\end{equation}
where $G_{\epsilon,N}$, $f_{\epsilon,N}$ and $\mathbf{S}_{\epsilon,N}$ are defined by
\begin{equation*}
G_{\epsilon,N} = g_{\epsilon}'(p_{\epsilon,N}) \, , \quad f_{\epsilon,N} = f_{\epsilon} - g_{\epsilon}(p_{\epsilon,N}) \, , \quad \mathbf{S}_{\epsilon,N} = \mathbf{S}_{\epsilon} - \nabla_{\mathbf{x}}p_{\epsilon,N} \, .
\end{equation*}
For each value of $N$, the problem (\ref{elliptic_non-linear_linearized_eps}) is of the same kind as (\ref{elliptic_quasi-linear_eps}). So we can solve it by applying the method described in the paragraph 2.1.2. \\
\indent This sequence of linearized problems can also be obtained from Newton's iterative method to solve
\begin{equation} \label{Newton_problem}
F_{\epsilon}(p_{\epsilon}) = 0 \, ,
\end{equation}
where the differential operator $F_{\epsilon}$ is defined as
\begin{equation*}
F_{\epsilon}(p) = -\nabla_{\mathbf{x}}\cdot \left(H_{\epsilon}\,(\mathbf{b}\otimes\mathbf{b})\,\cfrac{\nabla_{\mathbf{x}}p-\mathbf{S}_{\epsilon}}{\epsilon}\right) + g_{\epsilon}(p)-f_{\epsilon} \, .
\end{equation*}
Indeed, Newton's method for solving (\ref{Newton_problem}) writes
\begin{equation*}
DF_{\epsilon}(p_{\epsilon,N})(p_{\epsilon,N+1}-p_{\epsilon,N}) = -F_{\epsilon}(p_{\epsilon,N}) \, ,
\end{equation*}
where $DF_{\epsilon}(p)$ is the derivative in $p$ of the differential operator $F_{\epsilon}(p)$ and is of the form
\begin{equation*}
DF_{\epsilon}(p)(\delta) = -\nabla_{\mathbf{x}}\cdot \left(H_{\epsilon}\,(\mathbf{b} \otimes \mathbf{b}) \, \cfrac{\nabla_{\mathbf{x}}\delta}{\epsilon} \right) + g_{\epsilon}'(p)\times \delta  \, .
\end{equation*}

\section{Numerical method} \label{numet}
\setcounter{equation}{0}

In this section, we present a numerical method which allows to solve the diffusion problems (\ref{elliptic_quasi-linear_eps}) and (\ref{elliptic_non-linear_eps_intro}) by using the decomposition approaches we have presented. First, we introduce some notations which will be used for the construction of the scheme, then we present the scheme itself for the general linear case (\ref{elliptic_quasi-linear_eps}). Finally, we present the discretized version of Gummel's algorithm for the non-linear case.

\subsection{Notations and definitions}

We consider a uniform mesh $(x_{i},y_{j})$ defined by
\begin{equation*}
\begin{split}
x_{i} = x_{min} + i\,\Delta x \, , \quad y_{j} = y_{min}+j\,\Delta y \, , \;
\Delta x = \cfrac{x_{max}-x_{min}}{N_{x}+1} \, , \quad \Delta y = \cfrac{y_{max}-y_{min}}{N_{y}+1} \, ,
\end{split}
\end{equation*}
and we assume that the simulation domain is $\Omega = [x_{-1/2},x_{N_{x}+1/2}] \times [y_{-1/2},y_{N_{y}+1/2}]$. We also consider the following subsets of $\Z^{2}$:
\begin{equation*}
\begin{split}
I &= \{0,\dots,N_{x}\} \times \{0,\dots,N_{y}\} \, , \\
\overline{I} &= \{-1,\dots,N_{x}+1\} \times \{-1,\dots,N_{y}+1\} \, , \\
I_{*} &= \{0,\dots,N_{x}-1\} \times \{0,\dots,N_{y}-1\} \, , \\
\overline{I_{*}} &= \{-1,\dots,N_{x}\} \times \{-1,\dots,N_{y}\} \, ,
\end{split}
\end{equation*}
and we consider the notation
\begin{equation} \label{space_step}
h = \max(\Delta x, \Delta y) \, .
\end{equation}

\indent Since the decomposition method we have presented in paragraph 2.2 is based on variational formulations of the diffusion problem for $p_{\epsilon}$ and uses the duality between the operators $p \mapsto \mathbf{b} \cdot \nabla_{\mathbf{x}}p$ and $p \mapsto \nabla_{\mathbf{x}} \cdot (\mathbf{b} \, p)$, we choose to approach these differential operators by $\D_{h}$ and $\D_{h,*}$ respectively such that the duality property is preserved at the discrete level. For this purpose, we define $\D_{h}$ and $\D_{h,*}$ such that
\begin{equation} \label{def_Dh}
(\D_{h}\theta)_{i+1/2,j+1/2} = \mathbf{b}_{i+1/2,j+1/2} \cdot \left(
\begin{array}{c}
\cfrac{\theta_{i+1,j+1}-\theta_{i,j+1}+\theta_{i+1,j}-\theta_{i,j}}{2\Delta x} \\
\cfrac{\theta_{i+1,j+1}-\theta_{i+1,j}+\theta_{i,j+1}-\theta_{i,j}}{2\Delta y}
\end{array}
\right) \, ,
\end{equation}
for all $\theta = (\theta_{i,j})_{(i,j) \, \in \, \overline{I}}$, and
\begin{equation} \label{def_Dhstar}
\begin{split}
(\D_{h,*}\chi)_{i,j} &= \sum_{\alpha \in \{\pm 1\}} \cfrac{(b_{x}\chi)_{i+1/2,j+\alpha/2}-(b_{x}\chi)_{i-1/2,j+\alpha/2}}{2\Delta x} \\
&\qquad + \sum_{\alpha \in \{\pm 1\}} \cfrac{(b_{y}\chi)_{i+\alpha/2,j+1/2}-(b_{y}\chi)_{i+\alpha/2,j-1/2}}{2\Delta y} \, ,
\end{split}
\end{equation}
for all $\chi = (\chi_{i+1/2,j+1/2})_{(i,j) \, \in \, \overline{I_{*}}}$.

\subsection{Linear problems} \label{QL_FD}

We assume that the function $p \mapsto g_{\epsilon}(p)$ is given by
\begin{equation*}
g_{\epsilon}(p)(x,y) = G_{\epsilon}(x,y)\,p(x,y) \, ,
\end{equation*}
where $G_{\epsilon}$ is analytically known. We also assume that the functions $\mathbf{b}$, $H_{\epsilon}$, $f_{\epsilon}$ and $\mathbf{S}_{\epsilon}$ are analytically known and we consider the following notations:
\begin{equation*}
\begin{split}
f_{\epsilon,i,j} &= f_{\epsilon}(x_{i},y_{j}) \; , \\
G_{\epsilon,i,j} &= G_{\epsilon}(x_{i},y_{j}) \; , \\
G_{\epsilon,i+1/2,j+1/2} &= G_{\epsilon}(x_{i+1/2},y_{j+1/2}) \; , \\
H_{\epsilon,i+1/2,j+1/2} &= H_{\epsilon}(x_{i+1/2},y_{j+1/2}) \; , \\
\mathbf{S}_{\epsilon,i+1/2,j+1/2} &= \mathbf{S}_{\epsilon}(x_{i+1/2},y_{j+1/2}) \; , \\
\mathbf{b}_{i+1/2,j+1/2} &= \mathbf{b}(x_{i+1/2},y_{j+1/2}) \, .
\end{split}
\end{equation*}
Then, the diffusion problem (\ref{elliptic_quasi-linear_eps}) can be discretized under the following form:
\begin{equation} \label{elliptic_quasi-linear_FD}
\left\{
\begin{array}{ll}
\left(-\D_{h,*}(H_{\epsilon}\,(\D_{h}p_{\epsilon,app}-\mathbf{b}\cdot\mathbf{S}_{\epsilon}) + \epsilon\,G_{\epsilon}\,p_{\epsilon,app}\right)_{i,j} = \epsilon\,f_{\epsilon,i,j} \, , & \forall\,(i,j) \in I\,, \\
\left(H_{\epsilon} (\D_{h}p_{\epsilon,app} - \mathbf{b} \cdot \mathbf{S}_{\epsilon})\, (\mathbf{b}\cdot \bm{\nu}) \right)_{i+1/2,j+1/2} = 0 \, , & \forall\,(i,j) \in \overline{I_{*}} \backslash I_{*} \, ,
\end{array}
\right.
\end{equation}
and the approximation of $p_{\epsilon,app}$ of $p_{\epsilon}$ is computed at the points $(x_{i},y_{j}) \in \Omega$. \\
\indent Since $\D_{h}$ and $\D_{h,*}$ have been chosen to be dual operators, we follow the decomposition approach we have presented in Section 2.1.2 at a discrete level by using some discrete variational formulations of (\ref{elliptic_quasi-linear_FD}). Writing 
\begin{equation*}
p_{\epsilon,app,i,j} = \pi_{\epsilon,app,i,j} + q_{\epsilon,app,i,j} \, ,
\end{equation*}
with $\pi_{\epsilon,app}$ satisfying 
\begin{equation*}
(\D_{h}\pi_{\epsilon,app})_{i+1/2,j+1/2} = 0 \, , \quad \forall\, (i,j) \in I_{*} \, ,
\end{equation*}
$\pi_{\epsilon,app}$ and $q_{\epsilon,app}$ are completely defined by
\begin{equation}\label{elliptic_quasi-linear_FD:pi}
\pi_{\epsilon,app,i,j} = \cfrac{1}{G_{\epsilon,i,j}} \, \left[ f_{\epsilon,i,j} + \left(\D_{h,*}(G_{\epsilon}\,h_{\epsilon,app})\right)_{i,j}\right] \, , 
\end{equation}
\begin{equation*}
q_{\epsilon,app,i,j} = \cfrac{1}{G_{\epsilon,i,j}} \, \left(\D_{h,*}(G_{\epsilon}\,l_{\epsilon,app})\right)_{i,j} \, ,
\end{equation*}
where $h_{\epsilon,app} = (h_{\epsilon,app,i+1/2,j+1/2})_{(i,j) \, \in \, \overline{I_{*}}}$ and $l_{\epsilon,app} = (l_{\epsilon,app,i+1/2,j+1/2})_{(i,j) \, \in \, \overline{I_{*}}}$ are computed by inverting the following systems:
\begin{equation} \label{quasi-linear_h_FD}
\left\{
\begin{array}{ll}
-\left(\D_{h}\left(\cfrac{1}{G_{\epsilon}} \, \D_{h,*} (G_{\epsilon} \, h_{\epsilon,app}) \right)\right)_{i+1/2,j+1/2} & \\
\qquad \qquad \qquad \qquad \qquad = \left(\D_{h} \left( \cfrac{f_{\epsilon}}{G_{\epsilon}} \right)\right)_{i+1/2,j+1/2} \, , & \forall\,(i,j) \in I_{*} \, , \\
h_{\epsilon,app,i+1/2,j+1/2} = 0 \, , & \forall\,(i,j) \in \overline{I_{*}} \backslash I_{*} \, ,
\end{array}
\right.
\end{equation}
\begin{equation} \label{eq_L_quasi-linear_eps_FD}
\left\{
\begin{array}{ll}
\left(-\D_{h} \left( \cfrac{1}{G_{\epsilon}} \, \D_{h,*} (H_{\epsilon} \, L_{\epsilon,app}) \right)+ \epsilon \, L_{\epsilon,app}\right)_{i+1/2,j+1/2}  & \\
\qquad \qquad \qquad = - \epsilon\,\left(\D_{h} \left( \cfrac{f_{\epsilon}}{G_{\epsilon}} \right) - \mathbf{b}\cdot\mathbf{S}_{\epsilon}\right)_{i+1/2,j+1/2} \, , & \forall\,(i,j) \in I_{*}\,, \\
L_{\epsilon,app,i+1/2,j+1/2} = 0 \, , & \forall\,(i,j) \in \overline{I_{*}} \backslash I_{*} \, ,
\end{array}
\right.
\end{equation}
and
\begin{equation} \label{eq_l_quasi-linear_eps_FD}
\left\{
\begin{array}{ll}
-\left(\D_{h} \left( \cfrac{1}{G_{\epsilon}} \, \D_{h,*} (G_{\epsilon} \, l_{\epsilon,app}) \right)\right)_{i+1/2,j+1/2} & \\
\qquad \qquad \qquad \qquad = (L_{\epsilon,app} - \mathbf{b}\cdot\mathbf{S}_{\epsilon})_{i+1/2,j+1/2} \, , & \forall\,(i,j) \in I_{*}\,, \\
l_{\epsilon,app,i+1/2,j+1/2} = 0 \, , & \forall\,(i,j) \in \overline{I_{*}} \backslash I_{*} \, .
\end{array}
\right.
\end{equation}

\subsection{Non-linear problems}

In this paragraph, we detail the discretized version of Gummel's algorithm presented in Section 2.2. In order to initialize the loop, we compute the following initial datas:
\begin{equation*}
\begin{array}{rcll}
p_{\epsilon,0,i,j} &=& p_{\epsilon,0}(x_{i},y_{j}) \, , & \forall\,(i,j) \in \overline{I} \, , \\
\mathbf{S}_{\epsilon,i+1/2,j+1/2} &=& \mathbf{S}_{\epsilon}(x_{i+1/2},y_{j+1/2}) \, , & \forall\, (i,j) \in \overline{I_{*}} \, , \\
f_{\epsilon,i,j} &=& f_{\epsilon}(x_{i},y_{j}) \, , & \forall\,(i,j) \in I \, , \\
H_{\epsilon,i+1/2,j+1/2} &=& H_{\epsilon}(x_{i+1/2},y_{j+1/2}) \, , & \forall\,(i,j) \in \overline{I_{*}} \, . 
\end{array}
\end{equation*}

Then, the $N$-th iteration of Gummel's algorithm is set as follows:
\begin{itemize}
\item \underline{Step 1:} assuming that
\begin{equation*}
p_{\epsilon,N,i,j} \, , \qquad \forall\,(i,j) \in \overline{I} \, ,
\end{equation*}
are known, we have
\begin{equation*}
\begin{array}{rcll}
(\mathbf{b} \cdot \mathbf{S}_{\epsilon,N})_{i+1/2,j+1/2} &=& (\mathbf{b} \cdot \mathbf{S}_{\epsilon} - \D_{h}p_{\epsilon,N})_{i+1/2,j+1/2} \, , & \forall\,(i,j) \in \overline{I_{*}} \, , \\
f_{\epsilon,N,i,j} &=& f_{\epsilon,i,j} - g_{\epsilon}(p_{\epsilon,N,i,j}) \, , & \forall\,(i,j) \in I \, , \\
G_{\epsilon,N,i,j} &=& g_{\epsilon}'(p_{\epsilon,N,i,j}) \, , & \forall\,(i,j) \in I \, , \\
G_{\epsilon,N,i+1/2,j+1/2} &=& g_{\epsilon}'(p_{\epsilon,N,i+1/2,j+1/2}) \, , & \forall \, (i,j) \in \overline{I_{*}} \, ,
\end{array}
\end{equation*}
where $p_{\epsilon,N,i+1/2,j+1/2} = \cfrac{1}{4} \, (p_{\epsilon,N,i+1,j+1}+p_{\epsilon,N,i+1,j}+p_{\epsilon,N,i,j+1}+p_{\epsilon,N,i,j})$. \\

\item \underline{Step 2:} we compute $h_{\epsilon,N,i+1/2,j+1/2}$ and $l_{\epsilon,N,i+1/2,j+1/2}$ for all $(i,j) \in \overline{I_{*}}$ by solving
\begin{equation*} \label{non-linear_h_FD}
\left\{
\begin{array}{ll}
-\left(\D_{h}\left(\cfrac{1}{G_{\epsilon,N}} \, \D_{h,*} (G_{\epsilon,N} \, h_{\epsilon,N}) \right)\right)_{i+1/2,j+1/2} & \\
\qquad \qquad \qquad \qquad \qquad = \left(\D_{h} \left( \cfrac{f_{\epsilon,N}}{G_{\epsilon,N}} \right)\right)_{i+1/2,j+1/2} \, , & \forall\,(i,j) \in I_{*} \, , \\
h_{\epsilon,N,i+1/2,j+1/2} = 0 \, , & \forall\,(i,j) \in \overline{I_{*}} \backslash I_{*} \, ,
\end{array}
\right.
\end{equation*}
\begin{equation*} \label{eq_L_non-linear_eps_FD}
\left\{
\begin{array}{ll}
\left(-\D_{h} \left( \cfrac{1}{G_{\epsilon,N}} \, \D_{h,*} (H_{\epsilon} \, L_{\epsilon,N}) \right)+\epsilon \, L_{\epsilon,N}\right)_{i+1/2,j+1/2} & \\
\qquad \qquad \qquad = - \epsilon \left(\D_{h} \left( \cfrac{f_{\epsilon,N}}{G_{\epsilon,N}} \right) - \mathbf{b} \cdot \mathbf{S}_{\epsilon,N}\right)_{i+1/2,j+1/2} \, , & \forall\,(i,j) \in I_{*}\, , \\
L_{\epsilon,N,i+1/2,j+1/2} = 0 \, , & \forall\,(i,j) \in \overline{I_{*}} \backslash I_{*} \, ,
\end{array}
\right.
\end{equation*}
and
\begin{equation*} \label{eq_l_non-linear_eps_FD}
\left\{
\begin{array}{ll}
-\left(\D_{h} \left( \cfrac{1}{G_{\epsilon,N}} \, \D_{h,*} (G_{\epsilon,N} \, l_{\epsilon,N}) \right)\right)_{i+1/2,j+1/2} & \\
\qquad \qquad \qquad \qquad = (L_{\epsilon,N} - \mathbf{b}\cdot\mathbf{S}_{\epsilon,N})_{i+1/2,j+1/2} \, , & \forall\,(i,j) \in I_{*}\,, \\
l_{\epsilon,N,i+1/2,j+1/2} = 0 \, , & \forall\,(i,j) \in \overline{I_{*}} \backslash I_{*} \, .
\end{array}
\right.
\end{equation*}

\item \underline{Step 3:} we compute $\delta_{\epsilon,N,i,j}$ for all $(i,j) \in I$ by using
\begin{equation*}
\delta_{\epsilon,N,i,j} = \cfrac{1}{G_{\epsilon,N,i,j}} \, \left[ f_{\epsilon,N,i,j} + \left(\D_{h,*}\left(G_{\epsilon,N}\,(h_{\epsilon,N}+l_{\epsilon,N})\right) \right)_{i,j} \right] \, ,
\end{equation*}
and we obtain $p_{\epsilon,N+1,i,j}$ for all $(i,j) \in I$. \\

\item \underline{Step 4:} we compute $p_{\epsilon,N+1,i,j}$ for all $(i,j) \in \overline{I}\backslash I$ by using the boundary condition
\begin{equation*}
(\D_{h}p_{\epsilon,N+1})_{i+1/2,j+1/2}-(\mathbf{b}\cdot\mathbf{S}_{\epsilon})_{i+1/2,j+1/2} = 0 \, , \qquad \forall \, (i,j) \in \overline{I_{*}} \backslash I_{*} \, .
\end{equation*}

\end{itemize}

\section{Numerical investigations of the AP-scheme} \label{numres}

\setcounter{equation}{0}

This section is devoted to numerical investigations of the Asymptotic-Preserving scheme derived in Sections \ref{dec} and \ref{numet}. The validation procedure consists in manufacturing $p_\epsilon$, an  analytic solution of the model problem~\eqref{elliptic_non-linear_eps_intro} which is compared to the numerical approximation $p_{\epsilon,app}$ carried out thanks to  the AP-scheme. These experiments are performed in two dimensions using a uniform Cartesian mesh independent of the anisotropy direction. For simplicity purpose, the first numerical experiments are performed in the framework on the linear model, but the conclusions drawn from these investigations apply to the general non-linear problem.

\subsection{Numerical convergence of the scheme}\label{sec:convergence:linear}

\begin{figure}
   \centering
   \subfigure[$e_{p} = \left\|p_{\epsilon}-p_{\epsilon,app}\right\|_{\ell^{p}(I)} / \left\|p_{\epsilon}\right\|_{\ell^{p}(I)}$.\label{cv_inh_linear:error}]{\begin{minipage}[c]{0.49\textwidth}
      \centering \includegraphics[width=\textwidth]{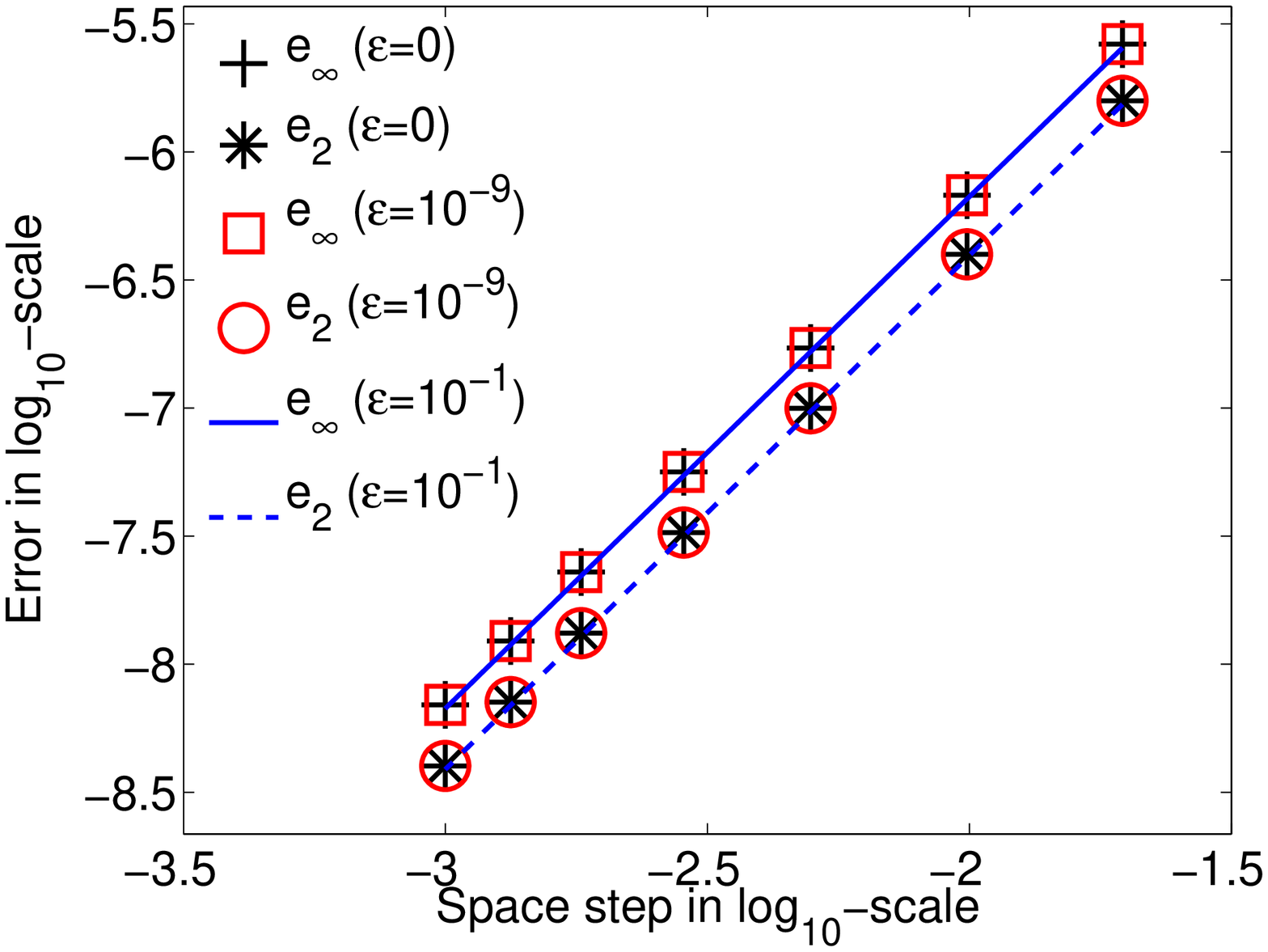}
     \end{minipage}}\hfill%
   \subfigure[$e_{p} = \left\|\D_{h}\pi_{\epsilon,app}\right\|_{\ell^{p}(I_{*})} / \left\|p_{\epsilon}\right\|_{\ell^{p}(I)}$. \label{cv_inh_linear:bdotgradpi}]{\begin{minipage}[c]{0.49\textwidth}\centering
       \includegraphics[width=\textwidth]{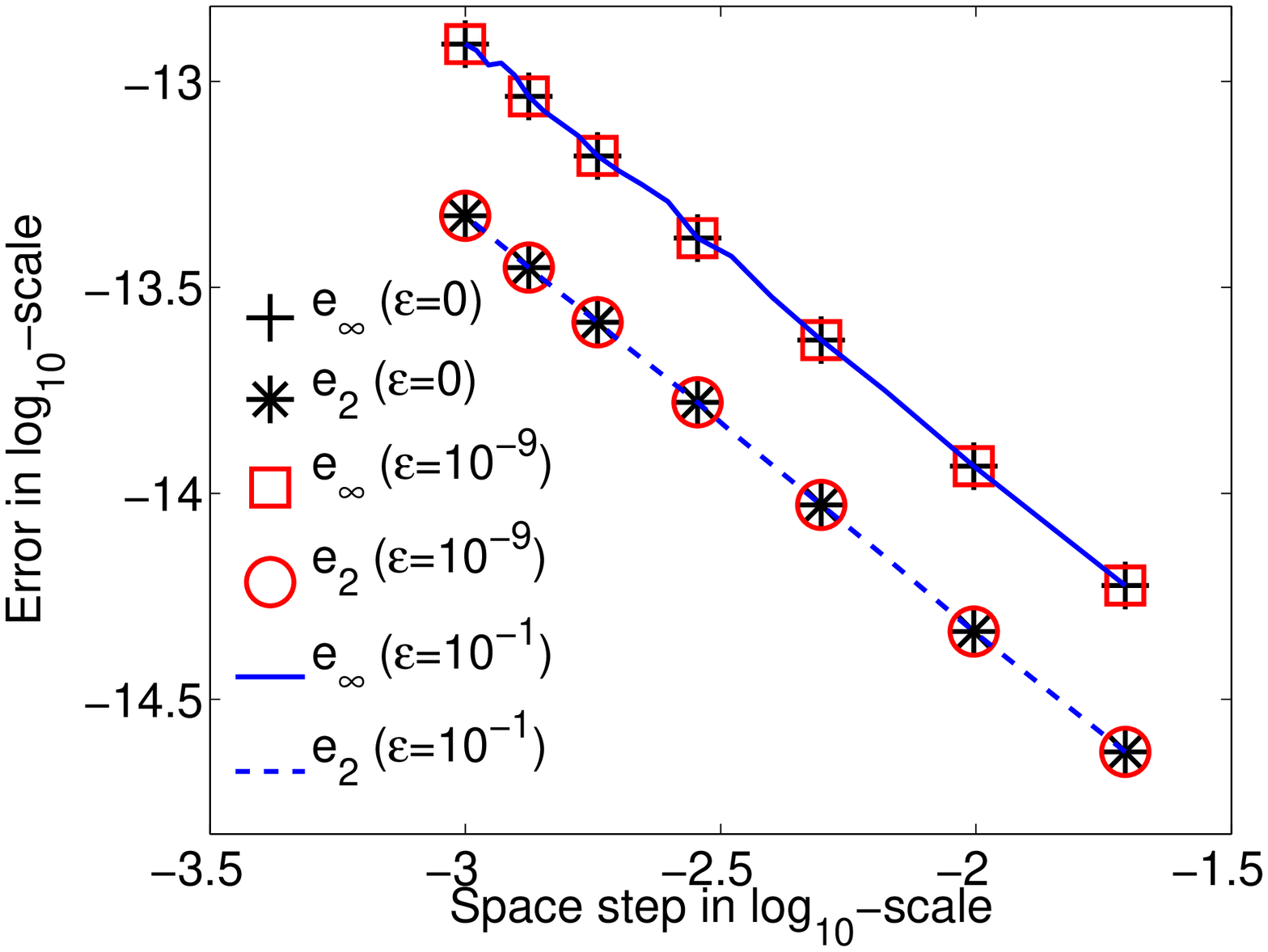}
     \end{minipage}}

  \caption{Relative error $\left\|p_{\epsilon}-p_{\epsilon,app}\right\|_{\ell^{p}(I)} / \left\|p_{\epsilon}\right\|_{\ell^{p}(I)}$ (left) and parallel gradient $\left\|\D_{h}\pi_{\epsilon,app}\right\|_{\ell^{p}(I_{*})} / \left\|p_{\epsilon}\right\|_{\ell^{p}(I)}$ (right) in $\ell^{2}$ and $\ell^{\infty}$ norms as functions of the space step $h$ for non-uniform data $G_{\epsilon}$, $H_{\epsilon}$ and $\mathbf{b}$ defined by (\ref{variable_GH})-(\ref{variable_b}), and different anisotropy strengths $\epsilon = 1$, $\epsilon = 10^{-9}$ and $\epsilon = 0$. \label{cv_inh_linear} }
%   \caption{$L^{2}$ and $L^{\infty}$ norms of the relative error $p_{\epsilon}-p_{\epsilon,app}$ (left) and of $\D_{h}\pi_{\epsilon,app}$ (right) as functions of the space step $h$: case with variable $G_{\epsilon}$, $H_{\epsilon}$ and $\mathbf{b}$ as defined in (\ref{variable_GH})-(\ref{variable_b}), and $\epsilon = 1$, $\epsilon = 10^{-9}$ and $\epsilon = 0$. \label{cv_inh_linear} }
\end{figure}

The first numerical tests aim at demonstrating the convergence of the AP-scheme regardless to the asymptotic parameter values. With this aim, an analytic solution is manufactured for the problem (\ref{elliptic_non-linear_eps_intro}) in the linear case, {\it i.e.} with $g_\epsilon(p) = G_\epsilon(\mathbf{x}) \, p(\mathbf{x})$. First, the expression for the anisotropy direction $\mathbf{b}$ and the functions $G_{\epsilon}$ and $H_{\epsilon}$ are defined on ${\Omega}=[1,2]\times[1,2]$ thanks to 
\begin{eqnarray} \label{variable_GH}
G_{\epsilon}(x,y) = H_{\epsilon}(x,y) = 1+\sin^{2}(x)\,\sin^{2}(y) \, ,
\\ \label{variable_b}
\mathbf{b} = (\sin\theta,-\cos\theta) \qquad \textnormal{with} \qquad \theta(x,y) = \arctan\left(\cfrac{y}{x}\right) \, .
\end{eqnarray}
Then the expression of $p_\epsilon$, as defined by
\begin{equation*}
p_{\epsilon}(x,y) = \cfrac{1}{1+x^{2}+y^{2}} \, ,
\end{equation*}
is used to analytically compute $f_{\epsilon}$ and $\mathbf{S}_{\epsilon}$ with
\begin{equation} \label{link_fSp}
f_{\epsilon} = G_{\epsilon}\,p_{\epsilon}\, , \qquad \mathbf{S}_{\epsilon} = \nabla_{\mathbf{x}}p_{\epsilon} \, .
\end{equation}
These definitions are inserted in the numerical method described in Section~\ref{QL_FD} to compute the numerical approximation $p_{\epsilon,app}$ finally compared to the exact solution $p_{\epsilon}$. The relative errors denoted $e_p$, $p \in \{2;\infty\}$, are defined by 
\begin{equation*}
e_{2} = \cfrac{\left\|p_{\epsilon}-p_{\epsilon,app}\right\|_{\ell^{2}(I)}}{\left\|p_{\epsilon}\right\|_{\ell^{2}(I)}} \, , \qquad e_{\infty} = \cfrac{\left\|p_{\epsilon}-p_{\epsilon,app}\right\|_{\ell^{\infty}(I)}}{\left\|p_{\epsilon}\right\|_{\ell^{\infty}(I)}} \, .
\end{equation*}
These quantities are displayed on Figure \ref{cv_inh_linear:error} as functions of the space step $h$ and for different anisotropy strengths $\epsilon = 10^{-1}$, $\epsilon = 10^{-9}$ and $\epsilon = 0$. A linear decrease of the errors is observed with the mesh refinement, the slope being equal to $2$, which is consistent with the definitions (\ref{def_Dh}) and (\ref{def_Dhstar}) of $\D_{h}$ and $\D_{h,*}$ as second order accurate approximations of the differential operators $\mathbf{b} \cdot \nabla_{\mathbf{x}}$ and $\nabla_{\mathbf{x}} \cdot (\mathbf{b}\, \cdot)$. Furthermore, this property holds for all considered values of $\epsilon$, including  $\epsilon = 0$. This demonstrates the $\epsilon$-invariance of the numerical scheme second order accuracy with respect  to the space step $h$.

\indent The ability of the scheme to compute a solution component $\pi_\epsilon$ with no gradient in the anisotropy direction is also investigated. The numerical approximation of $\pi_\epsilon$, $\pi_{\epsilon,app}$ provided by (\ref{elliptic_quasi-linear_FD})-(\ref{elliptic_quasi-linear_FD:pi}), should verify a discrete analogous of the property $\mathbf{b}\,\cdot \nabla_{\mathbf{x}} \pi_\epsilon =0$. This is analyzed thanks to Figure~\ref{cv_inh_linear:bdotgradpi}, where the evolution of $\left\|\D_{h}\pi_{\epsilon,app}\right\|_{\ell^{p}(I_{*})} / \left\|p_{\epsilon}\right\|_{\ell^{p}(I)}$ as a function of the space step is displayed for $p=2,\infty$,  $\epsilon = 10^{-1}$, $\epsilon = 10^{-9}$ and $\epsilon = 0$. Note that the quantity $\D_{h}\pi_{\epsilon,app}$  is the residual of the linear system solved to compute the solution of  \eqref{elliptic_quasi-linear_FD}, and consequently characterizes the precision of the linear system solver. For these test cases, a sparse direct solver being used \cite{mumps}, the accuracy is very close to the computer arithmetic precision, at least for small linear system sizes. This precision is observed to  deteriorate moderately with the increase of the system size which explains the growth of the error with vanishing mesh sizes. However this does not affect the precision of the scheme, as demonstrated by the results of Figure~\ref{cv_inh_linear:error}.

\subsection{Anisotropy angle influence on the method accuracy}

In this section, we quantify the sensitivity of the numerical method with respect to the anisotropy direction variations. More precisely, we wish to analyze the accuracy of the method as a function of $\alpha$, the angle measured between the anisotropy direction and the first direction (associated to the first coordinate). The anisotropy direction $\mathbf{b}$ is assumed to be uniform and defined as
\begin{equation*}
\mathbf{b} = (\sin\alpha,-\cos\alpha) \, , \hspace*{5 mm} \alpha \in [0,\frac{\pi}{2}].
\end{equation*}
\indent In order to manufacture an analytic solution for the problem, we introduce a system of coordinates which is adapted to $\mathbf{b}$. These coordinates are denoted $(X,Y)$ and are deduced from $(x,y)$ by the relations
\begin{equation}\label{eq:def:chgt:var}
X = x\,\cos\alpha + y\,\sin\alpha \, , \qquad Y = x\,\sin\alpha-y\,\cos\alpha \, .
\end{equation}
In these coordinates, the linear diffusion problem (\ref{elliptic_quasi-linear_eps}) writes
\begin{equation} \label{eq_quasi-linear_eps_aligned}
\left\{
\begin{array}{ll}
-\D_{Y}(\mathcal{H}_{\epsilon}\,\D_{Y}\mathcal{P}_{\epsilon}) + \epsilon\,\mathcal{G}_{\epsilon}\,\mathcal{P}_{\epsilon} = \epsilon\,\mathcal{F}_{\epsilon} - \D_{Y}(\mathcal{H}_{\epsilon}\,\mathcal{B}\cdot\mathcal{S}_{\epsilon}) \, , & \textnormal{in $\Omega$,} \\
\D_{Y}\mathcal{P}_{\epsilon} = \mathcal{B}\cdot\mathcal{S}_{\epsilon} \, , &\textnormal{on $\D\Omega$,}
\end{array}
\right.
\end{equation}
with $\mathcal{P}_{\epsilon}(X,Y) = p_{\epsilon}(x,y)$, $\mathcal{H}_{\epsilon}(X,Y) = H_{\epsilon}(x,y)$, $\mathcal{G}_{\epsilon}(X,Y) = G_{\epsilon}(x,y)$, $\mathcal{F}_{\epsilon}(X,Y) = f_{\epsilon}(x,y)$, $\mathcal{S}_{\epsilon}(X,Y) = \mathbf{S}_{\epsilon}(x,y)$ and $\mathcal{B}(X,Y) = \mathbf{b}(x,y)$. It is straightforward to verify that the function $\mathcal{P}_{\epsilon}$ given by
\begin{equation*}
\mathcal{P}_{\epsilon}(X,Y) = \sin(X) + \cfrac{1}{\mathcal{G}_{\epsilon}(X,Y)} \, \D_{Y}(\mathcal{G}_{\epsilon}\,\mathcal{L}_{\epsilon}) \, ,
\end{equation*}
is the solution of (\ref{eq_quasi-linear_eps_aligned}) provided that $\mathcal{F}_{\epsilon}$ and $\mathcal{S}_{\epsilon}$ satisfy
\begin{equation*}
\mathcal{F}_{\epsilon} = \mathcal{G}_{\epsilon}\,\mathcal{P}_{\epsilon} \, , \qquad \mathcal{B} \cdot \mathcal{S}_{\epsilon} = \D_{Y} \mathcal{P}_{\epsilon} \, ,
\end{equation*}
and where $\mathcal{L}_{\epsilon}(X,Y) = l_{\epsilon}(x,y)$ with ${l_{\epsilon}}_{|_{\D\Omega}} = 0$. This requirement is met by the following definition 
\begin{equation*}
l_{\epsilon}(x,y) = \sin\left(\cfrac{2\pi\,(x-x_{-1/2})}{x_{N_{x}+1/2}-x_{-1/2}} \right) \, \sin\left(\cfrac{2\pi\,(y-y_{-1/2})}{y_{N_{y}+1/2}-y_{-1/2}} \right) \, ,
\end{equation*}
which ensures $l_{\epsilon}(x_{i+1/2},y_{j+1/2}) = 0$ for any $(i,j) \in \overline{I_{*}} \backslash I_{*}$. 
The problem is stated in Cartesian coordinates thanks to the change of variables~\eqref{eq:def:chgt:var} yielding to $p_{\epsilon}(x,y) = \pi_{\epsilon}(x,y) + q_{\epsilon}(x,y)$
with
\begin{eqnarray*}
\pi_{\epsilon}(x,y) = \sin\left(x\,\cos\alpha + y\,\sin\alpha\right) \, ,
\hspace*{2 mm}
q_{\epsilon}(x,y) = \cfrac{1}{G_{\epsilon}(x,y)} \, \nabla_{\mathbf{x}} \cdot \left(G_{\epsilon}(x,y)\,\mathbf{b}(x,y)\,l_{\epsilon}(x,y) \right) \, ,
\end{eqnarray*}
the other coefficients being manufactured similarly with $G_{\epsilon}$ and $H_{\epsilon}$ given by (\ref{variable_GH}).

\begin{figure}
  \centering
  \subfigure[$\left\|p_{\epsilon}-p_{\epsilon,app}\right\|_{\ell^{p}(I)} / \left\|p_{\epsilon}\right\|_{\ell^{p}(I)}$, $\epsilon = 10^{-3}$. \label{angle_GH_variable:error_10m3}]{\begin{minipage}[c]{0.49\textwidth}\centering
     \includegraphics[width=\textwidth]{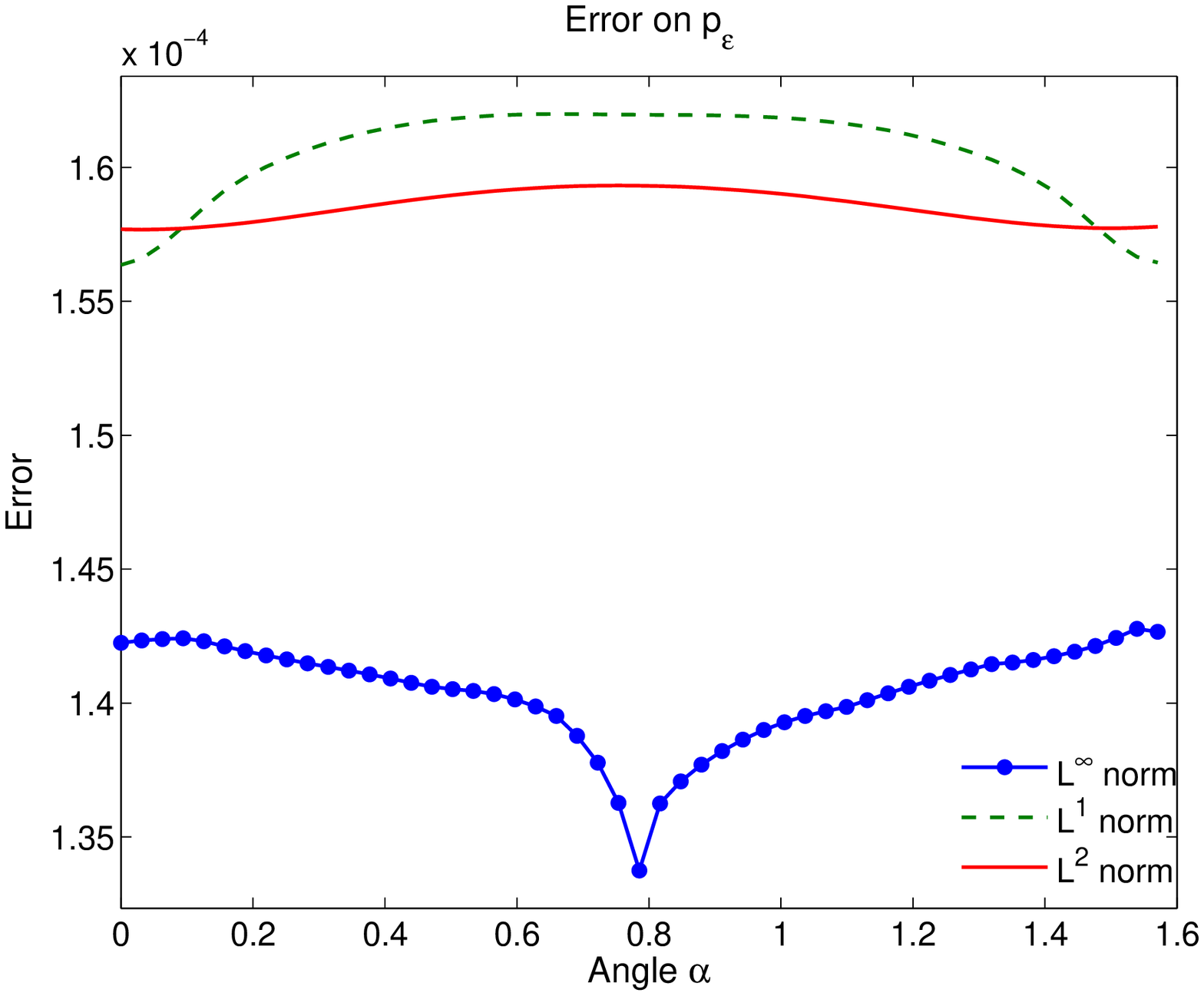}
    \end{minipage}}\hfill%
  \subfigure[$\left\|p_{\epsilon}-p_{\epsilon,app}\right\|_{\ell^{p}(I)} / \left\|p_{\epsilon}\right\|_{\ell^{p}(I)}$, $\epsilon = 10^{-8}$. \label{angle_GH_variable:error_10m8}]{\begin{minipage}[c]{0.49\textwidth}\centering
      \includegraphics[width=\textwidth]{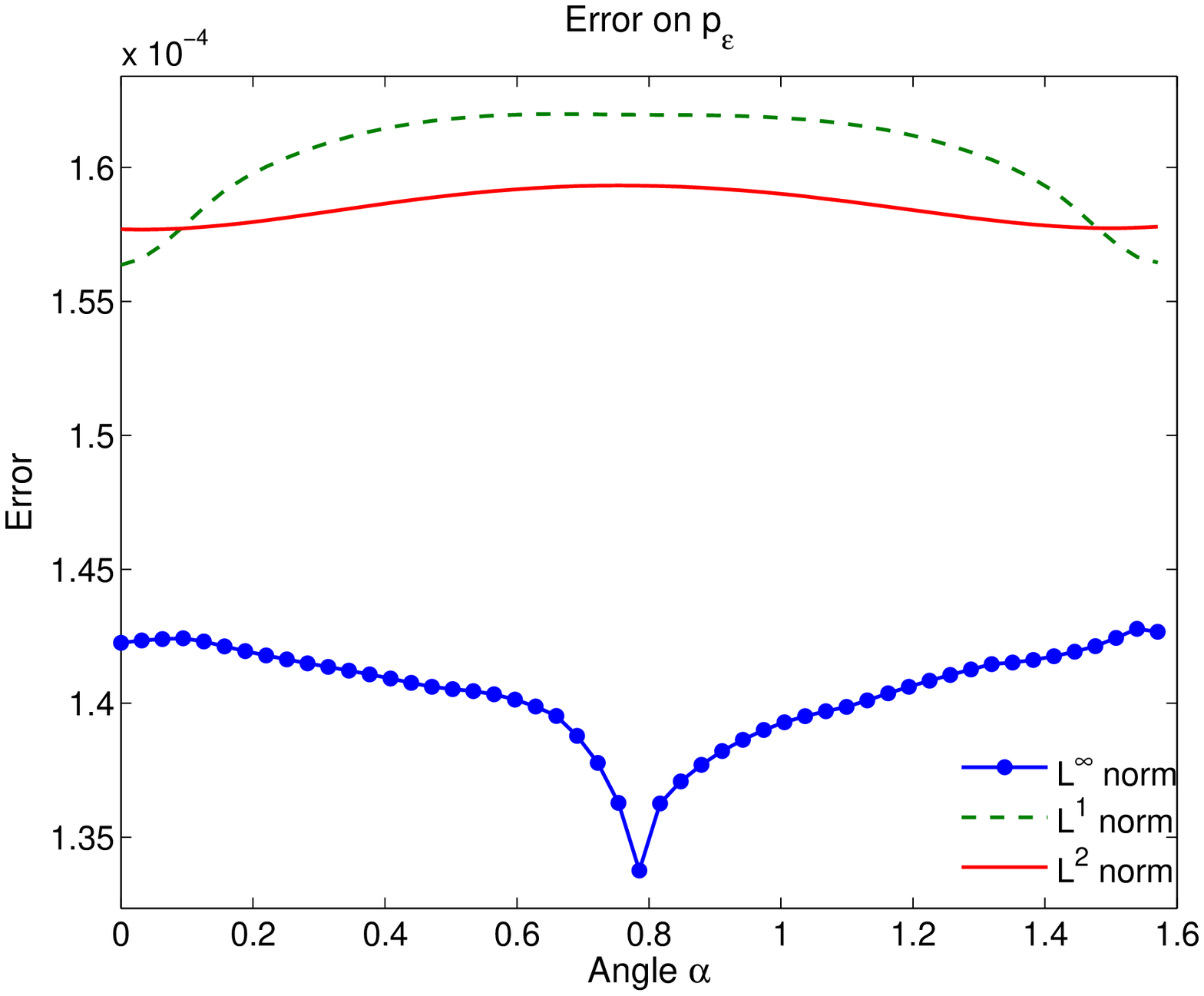}
    \end{minipage}}
  \caption{Relative error $\left\|p_{\epsilon}-p_{\epsilon,app}\right\|_{\ell^{p}(I)} / \left\|p_{\epsilon}\right\|_{\ell^{p}(I)}$ ($p=1,2,\infty$) as a function of $\alpha$: case with $G_{\epsilon}(x,y) = H_{\epsilon}(x,y) = 1+\sin^{2}(x)\,\sin^{2}(y)$ and $\epsilon = 10^{-3}$ (left) and $\epsilon = 10^{-8}$ (right). \label{angle_GH_variable}}
\end{figure}

In the following tests, the computation domain $\Omega = [1,2] \times [1,2]$ is discretized thanks to  a uniform mesh constituted of $200 \times 200$ cells. The relative approximation error as a function of the angle $\alpha$ is displayed on Figure \ref{angle_GH_variable} for different norms. The numerical method accuracy is observed to remain almost unaltered by the anisotropy direction changes. More precisely, we observe a variation of the relative errors in norms $\ell^{1}$ and $\ell^{2}$ lower than $4\%$ and a variation of $\ell^{\infty}$ lower than $7\%$. Furthermore, these observations are redundant for several values of $\epsilon$: in Figure \ref{angle_GH_variable}, we have considered $\epsilon = 10^{-3}$ and $\epsilon = 10^{-8}$ and the obtained error curves are very close. Note that other experiments have been carried out for anisotropy strengths ranging from $\epsilon=0$ to $\epsilon=1$ and with other definitions of $G_{\epsilon}$ and of $H_{\epsilon}$, with comparable results. The curves being very similar to that of Figure \ref{angle_GH_variable}, these plots are omitted.

\indent These observations confirm one of the main ideas of the present paper: the accuracy of the method is almost independent of the anisotropy direction relatively to the grid, \textit{i.e.} the mesh over $\Omega$ can be constructed whatever the anisotropy direction and strength, without a significant loss of accuracy.

\subsection{Convergence of Gummel's loop}

\indent The third test sequence is devoted to the convergence of the linear problems sequence defined in Sections 2.2 and 3.3 for solving the non-linear model (\ref{elliptic_non-linear_eps_intro}). The process detailed in the preceding sections is again implemented to manufacture an analytic solution for the non-linear problem.

The computational domain remains $\Omega = [1,2]\times[1,2]$, the anisotropy direction $\mathbf{b}$ is a function of the space variables whose expression is given by equation \eqref{variable_b},  $H_{\epsilon}$ and $g_{\epsilon}$ being defined as
\begin{equation} \label{def_Hepsgeps}
H_{\epsilon}(x,y) = 1+\cos^{2}(x)\,\cos^{2}(y) \, , \qquad g_{\epsilon}(p) = p^{6} \, .
\end{equation}
Note that this choice of $g_{\epsilon}(p)$ introduces a severe non-linearity in the problem. Several tests have also been performed with other definitions of $g_{\epsilon}(p)$, for instance
\begin{displaymath}
g_{\epsilon}(p) = -p\,(1-p)\,\left(p-\frac{1}{2}\right)\, ,
\end{displaymath}
which defines an anisotropic diffusion-reaction equation similar to the steady-state Allen-Cahn equation (see \cite{Benes-Yazaki-Kimura,Elliot-Schatzle}) used in phase transition problems. These tests produce results almost identical to the results which are obtained when $g_{\epsilon}(p) = p^{6}$ so we only consider the strongly non-linear reaction term defined in (\ref{def_Hepsgeps}) within the presentation of the numerical results in the next lines. \\
\indent The solution $p_{\epsilon}$ is constructed thanks to a cubic spline $S$, precisely 
\begin{equation} \label{def_peps_spline}
p_{\epsilon}(x,y) = 1+S\left(\cfrac{x-x_{mid}}{L_{x}}\right) \, S\left(\cfrac{y-y_{mid}}{L_{y}}\right) \, ,
\end{equation}
with $S(z) = 0$ for $|z| \notin [0,2]$  and 
\begin{equation} \label{def_spline}
S(z) = \left\{
\begin{array}{ll}
\frac{1}{6} \, (2-|z|^{3}) \, , & \textnormal{if $1\leq |z| \leq 2$,} \\[0.3em]
\frac{2}{3} - |z|^{2} + \frac{1}{2}\,|z|^{3}\, , \quad & \textnormal{if $0 \leq |z| < 1$,} \\ 
\end{array}
\right.
\end{equation}
with $(x_{mid},y_{mid}) = (\frac{3}{2},\frac{3}{2})$ and $L_{x} = L_{y} = {1/10}$.
To analyze the convergence  with respect to the number of Gummel's iterations, the sequence is initiated with $p_{\epsilon,0}$, a perturbation of the non-linear problem solution, reading
\begin{equation} \label{def_peps0}
p_{\epsilon,0}(x,y) = p_{\epsilon}(x,y) + \eta \, \max\left(0,1-\mu\,(x-x_{mid})^{2}-\mu\,(y-y_{mid})^{2}\right) \, ,
\end{equation}
where $\mu$ and $\eta\,\mu$ are parameters controlling the support and the magnitude of the perturbation. Since Gummel's method is constructed on a linearization of the problem its convergence cannot be guaranteed with a poor estimation of the solution as initial guess. It means that the parameters $\mu$ and $\eta$ cannot be chosen completely arbitrarily: indeed, several simulations have been performed, all with the same parameters except $\eta$ ranging in $\{0,10,20, 40, 60, 100, 1000\}$ and $\mu$ ranging in $\{1,10,60,100\}$ and it has been observed that Gummel's method does not converge as $N \to \infty$ when $\eta$ is larger than $10^{2}$. Concerning the parameter $\mu$, the simulation sequence reveals that the convergence of Gummel's method is almost not affected by the amplitude of $\mu$.

\begin{figure}
 \centering
 \subfigure[$E_{2,N}$, $M_{100}$ and $M_{1000}$ meshes. \label{convergence_inN_figs:error}]{\begin{minipage}[c]{0.49\textwidth}
     \includegraphics[width=\textwidth]{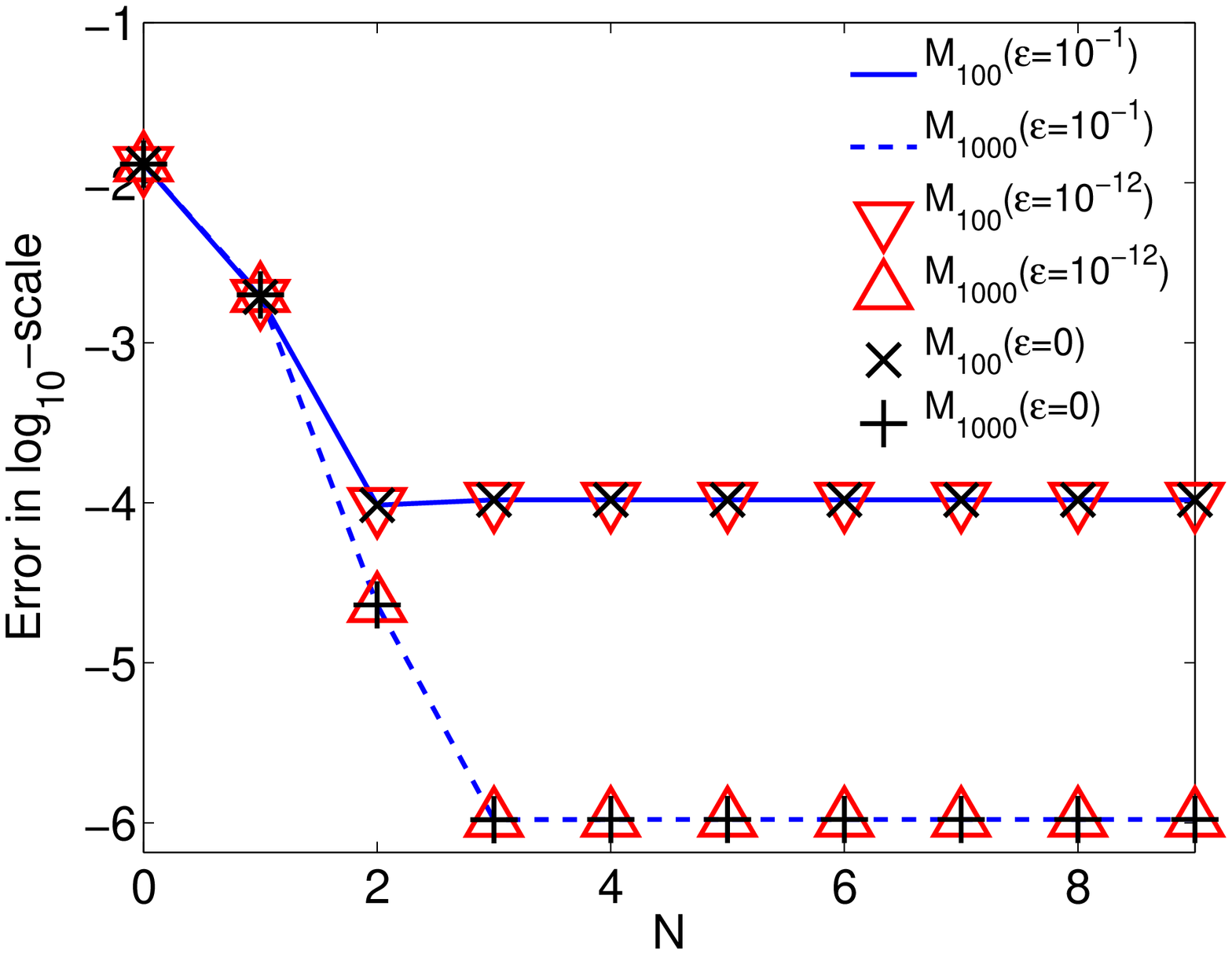}
   \end{minipage}}\hfill%
 \subfigure[$\|\delta_{\epsilon,N,app}\|_{\ell^{2}(I)} / \|p_{\epsilon}\|_{\ell^{2}(I)}$, $M_{100}$ mesh. \label{convergence_inN_figs:corrector}]{\begin{minipage}[c]{0.49\textwidth}
     \includegraphics[width=\textwidth]{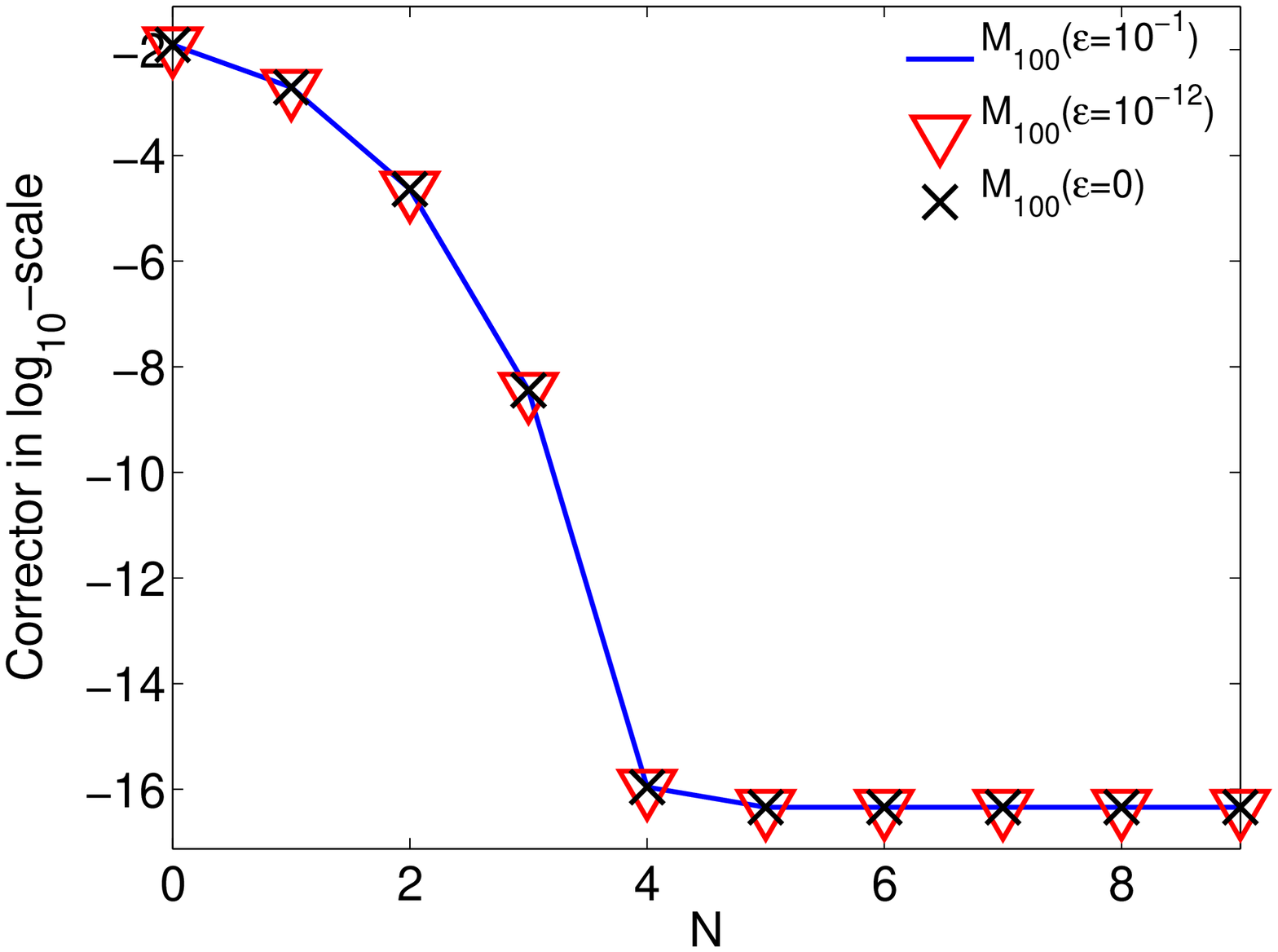}
   \end{minipage}}
 \caption{Gummel's iteration convergence: evolution in $\log_{10}$-scales of the relative error $E_{2,N} = \|p_{\epsilon}-p_{\epsilon,N_{f},app}\|_{\ell^{2}(I)} / \|p_{\epsilon}\|_{\ell^{2}(I)}$ (left) and of $\|\delta_{\epsilon,N,app}\|_{\ell^{2}(I)} / \|p_{\epsilon}\|_{\ell^{2}(I)}$ (right) as functions of th iteration number $N$ for the non-linear problem. The simulations are performed with the uniform meshes $M_{100}$ and $M_{1000}$ composed of $100\times100$ and $1000\times 1000$ cells (for the corrector norm, the values being very similar for both meshes, only those of the coarsest one are displayed).} \label{convergence_inN_figs}
\end{figure}

The successive relative errors measured between the iterates of the Gummel's loop and the exact solution are plotted on Figure~\ref{convergence_inN_figs:error}. The computations are carried out on two different meshes,  $M_{100}$ and $M_{1000}$ with $100\times100$ and $1000\times 1000$ cells, with $\eta = 0.1$ and $\mu = 60$ and for anisotropy strengths including $\epsilon=0$.  Along with the graphical representation of the solution approximation error, the evolution of the corrector norm relative to that of the solution, namely the quantity $\|\delta_{\epsilon,N,app}\|_{\ell^{2}(I)} / \|p_{\epsilon}\|_{\ell^{2}(I)}$, is also plotted in Figure \ref{convergence_inN_figs:corrector}. These last results being almost identical for both meshes, the plot related to the finest mesh is omitted in this figure.

\begin{table}
\begin{center}

\begin{tabular}{|c|c|c|c|c|c|}
\hline
$\epsilon$ & Error & $M_{100}$ & $M_{200}$ & $M_{500}$ & $M_{1000}$ \\
\hline
& $E_{1,N_{f}}$ & $3.9452 \times 10^{-5}$ & $9.8116 \times 10^{-6}$ & $1.5673 \times 10^{-6}$ & $3.9166 \times 10^{-7}$ \\
\cline{2-6} 
$10^{-1}$ & $E_{2,N_{f}}$ & $1.0446 \times 10^{-4}$ & $2.6188 \times 10^{-5}$ & $4.1988 \times 10^{-6}$ & $1.0505 \times 10^{-6}$ \\
\cline{2-6} 
& $E_{\infty,N_{f}}$ & $6.0730 \times 10^{-4}$ & $1.5793 \times 10^{-4}$ & $2.5942 \times 10^{-5}$ & $6.5451 \times 10^{-6}$ \\
\hline
& $E_{1,N_{f}}$ & $3.9796 \times 10^{-5}$ & $9.8969 \times 10^{-6}$ & $1.5808 \times 10^{-6}$ & $3.9504 \times 10^{-7}$ \\
\cline{2-6} 
$10^{-12}$ & $E_{2,N_{f}}$ & $1.0496 \times 10^{-4}$ & $2.6311 \times 10^{-5}$ & $4.2184 \times 10^{-6}$ & $1.0554 \times 10^{-6}$ \\
\cline{2-6} 
& $E_{\infty,N_{f}}$ & $6.1098 \times 10^{-4}$ & $1.5885 \times 10^{-4}$ & $2.6087 \times 10^{-5}$ & $6.5815 \times 10^{-6}$ \\
\hline
& $E_{1,N_{f}}$ & $3.9796 \times 10^{-5}$ & $9.8969 \times 10^{-6}$ & $1.5808 \times 10^{-6}$ & $3.9504 \times 10^{-7}$ \\
\cline{2-6} 
$0$ & $E_{2,N_{f}}$ & $1.0496 \times 10^{-4}$ & $2.6311 \times 10^{-5}$ & $4.2184 \times 10^{-6}$ & $1.0554 \times 10^{-6}$ \\
\cline{2-6} 
& $E_{\infty,N_{f}}$ & $6.1098 \times 10^{-4}$ & $1.5885 \times 10^{-4}$ & $2.6087 \times 10^{-5}$ & $6.5815 \times 10^{-6}$ \\
\hline
\end{tabular}
\end{center}
\caption{Relative error $E_{p,N_{f}} = \|p_{\epsilon}-p_{\epsilon,N_{f},app}\|_{\ell^{p}(I)} / \|p_{\epsilon}\|_{\ell^{p}(I)}$ ($p = 1,2,\infty$) for the non-linear problem defined by $g_{\epsilon}(p) = p^{6}$. The computations are carried out on uniform meshes $M_{k}$ constituted of $k \times k$ cells ($k = 100, 200, 500, 1000$) with several values of $\epsilon$ and after a number of iteration of Gummel's loop $N_f$ large enough for the convergence to be effective.} \label{convergence_inN_tab}
\end{table}

\indent In spite of the large perturbation amplitude,  Gummel's iterative method converges in a small number of iterations, for both meshes and for all $\epsilon$-values. The corrector term $\delta_{\epsilon,N,app}$ rapidly decreases to reach the computer precision threshold ($10^{-15}$) after 4 iterations. In the same time, the relative error also decreases but the approximation is not improved by subsequent iterations, the error remaining constant for iteration numbers greater than 4. At this stage, the precision of the approximation is not limited by the linearization process of the Gummel's loop anymore, but by the discretization error of the linearized problem, explaining the plateau described by the error. To document this analyzis further, we summarize in Table \ref{convergence_inN_tab} the values of the relative error measured between the exact solution and the approximation obtained after $N_f$ iterations of the Gummel's loop. This quantity is referred to as $E_{p,N_{f}}$ ($p=1,2,\infty$) and computed for $N_{f}$  large enough to ensure that the plateau above mentioned is reached. For the investigations carried out, this requirement is met as soon as $N_f\geq 4$. The approximation error  $E_{p,N_f}$ is observed to quadratically decrease with the space mesh: the error norms related to the computations performed on a $1000\times1000$ mesh  are for instance $10^2$ times as small as those carried out on a mesh with $100\times 100$ cells. This is a consequence of the second order accurate discretization of the spatial operator already outlined in section~\ref{sec:convergence:linear}. Finally, the results of Table~\ref{convergence_inN_tab} also demonstrate the independence of the numerical method precision with respect to the anisotropy intensity.

\subsection{Highlight of the scheme Asymptotic-Preserving property}

\begin{figure}
 \centering
 \subfigure[$E_{\epsilon}$, $M_{200}$ and $M_{1000}$ meshes. \label{convergence_ineps_nonlinear_figs:exact}]{\begin{minipage}[c]{0.49\textwidth}
     \includegraphics[width=\textwidth]{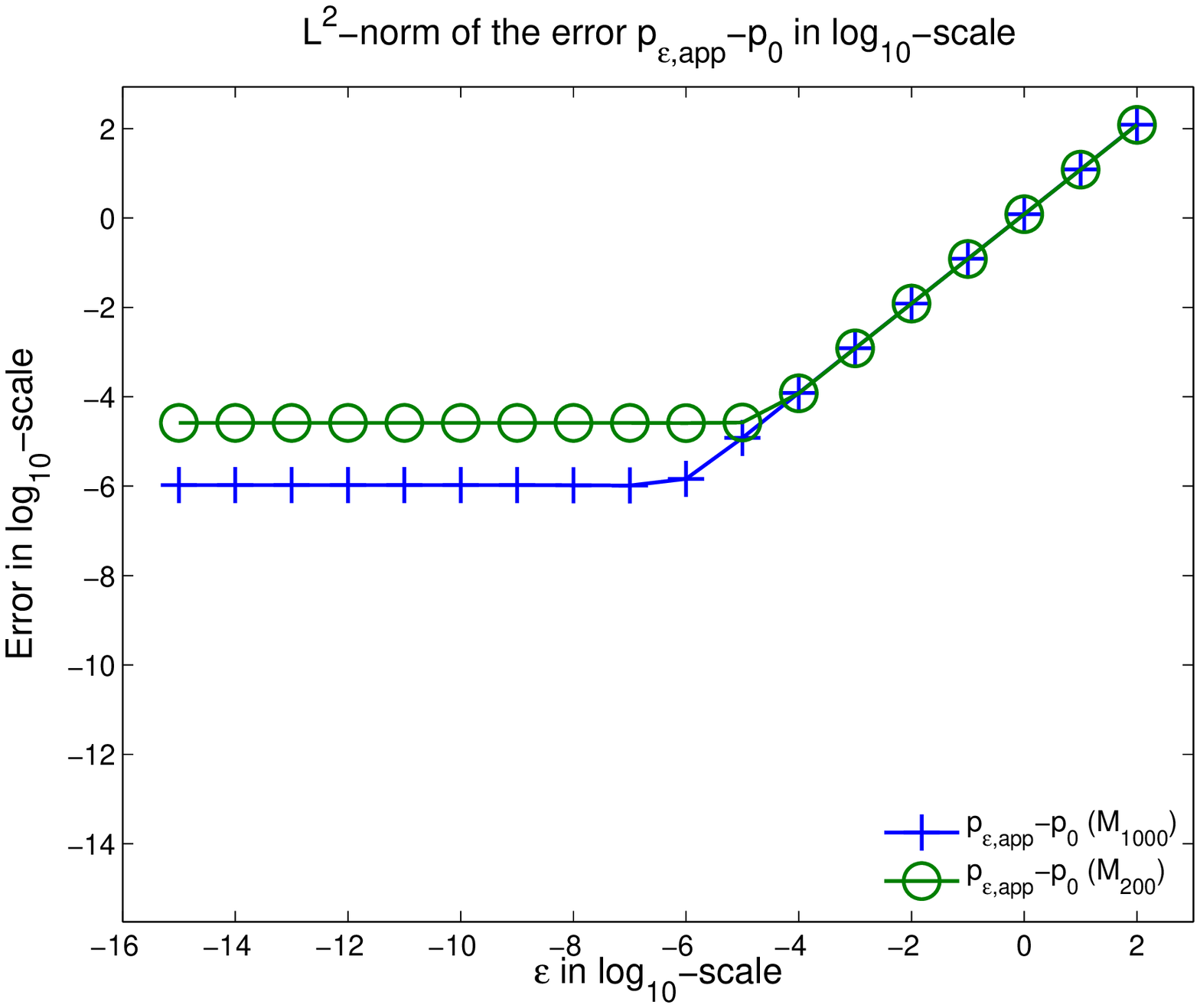}
   \end{minipage}}\hfill%
 \subfigure[$E_{\epsilon,app}$, $M_{200}$ and $M_{1000}$ meshes. \label{convergence_ineps_nonlinear_figs:approx}]{\begin{minipage}[c]{0.49\textwidth}
     \includegraphics[width=\textwidth]{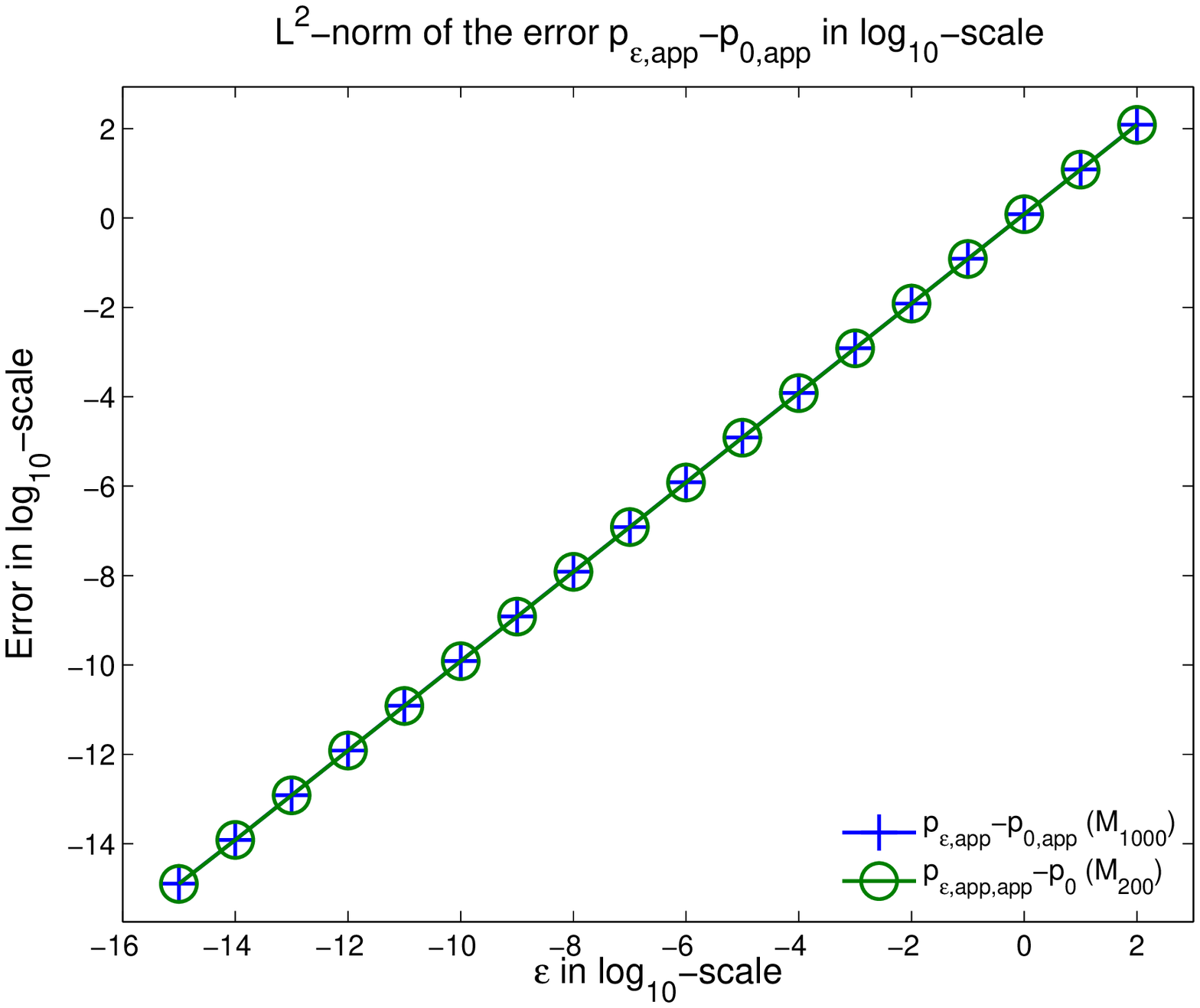}
   \end{minipage}}

  \caption{Evolution in $\log_{10}$-scales of $E_{\epsilon} = \left\|p_{\epsilon,app}-p_{0}\right\|_{\ell^{2}(I)} / \left\|p_{0}\right\|_{\ell^{2}(I)}$ (left) and of $E_{\epsilon,app} = \left\|p_{\epsilon,app}-p_{0,app}\right\|_{\ell^{2}(I)} / \left\|p_{0}\right\|_{\ell^{2}(I)}$ (right) as functions of $\epsilon$ computed on uniform meshes constituted of $200 \times 200$ cells ($M_{200}$) and $1000 \times 1000$ cells ($M_{1000}$) are considered. The simulations are performed with $\mathbf{b}$, $H_{\epsilon}$, $g_{\epsilon}$ defined by (\ref{def_Hepsgeps}) and (\ref{variable_b}).} \label{convergence_ineps_nonlinear_figs}
\end{figure}

These last experiments are devoted to illustrate the Asymptotic-Preserving property of the numerical method,  \textit{i.e.} its ability to compute an accurate approximation of $p_{0}$, the solution of the limit problem~\eqref{eq:def:limit:problem}. The solution of the problem is constructed as  a sequence $(p_{\epsilon})_{\epsilon\,>\,0}$ defined by
\begin{equation*}
p_{\epsilon} = p_{0} + \epsilon\,\tilde{p}_{\epsilon}^1 \, ,
\end{equation*}
with 
\begin{align}
p_{0}(x,y) &= 1+S\left(\cfrac{x-x_{mid}}{L_{x}}\right) \, S\left(\cfrac{y-y_{mid}}{L_{y}}\right) \, ,\\
\tilde{p}_{\epsilon}^1(x,y) &= \max\left(0, \cos\left(\cfrac{2\pi\,(x-x_{mid})}{L_{x}}\right)\, \cos\left(\cfrac{2\pi\,(y-y_{mid})}{L_{y}}\right) \right) \, .
\end{align}
The functions $g_{\epsilon}$, $H_{\epsilon}$ and $\mathbf{b}$ are defined as in the previous test sequence, the initial guess for Gummel's loop being constructed following (\ref{def_peps0}) using the same perturbation. We now wish to evaluate the error measured between the exact solution of the limit problem $p_0$ and the approximation computed thanks to the AP-scheme for vanishing $\epsilon$. This error, denoted $E_{\epsilon}$ and defined as 
\begin{equation*}
E_{\epsilon} = {\left\|p_{\epsilon,app}-p_{0}\right\|_{\ell^{2}(I)}}/{\left\|p_{0}\right\|_{\ell^{2}(I)}}  \, ,
\end{equation*}
is plotted on Figure~\ref{convergence_ineps_nonlinear_figs:exact} as a function of $\epsilon$. The data represented on this figure are obtained after convergence of the Gummel's loop. Two regimes can be identified. The first one is related to the largest values of $\epsilon$ for which a linear decrease of the error is observed. The second one is a plateau whose value depends on the mesh step $h$, this value being lower for refined meshes. Precisely we note a quadratic decrease of this value with the mesh size. To explain these features, we use the following identity
\begin{equation*}
p_{\epsilon,app}-p_{0} = p_{\epsilon,app}-p_{0,app} + {p_{0,app}-p_{0}} \, . %_{= \mathcal{O}(h^{2})} \, .
\end{equation*}
This yields $ E_{\epsilon} \leq E_{\epsilon,app} + e_0$ where $e_0 = \left\|p_{0,app}-p_{0}\right\|_{\ell^{2}(I)} / \left\|p_{0}\right\|_{\ell^{2}(I)}$ represents the approximation error of $p_0$, $p_{0,app}$ being the numerical approximation of $p_{0}$ provided by the AP-scheme with $\epsilon = 0$, and $E_{\epsilon,app} = \left\|p_{\epsilon,app}-p_{0,app}\right\|_{\ell^{2}(I)} / \left\|p_{0}\right\|_{\ell^{2}(I)}$. 
The error $E_{\epsilon}$ linearly decreases with $\epsilon$ as long as the approximation error $e_0$ is negligible compared to $E_{\epsilon,app}$ (see Figure~\ref{convergence_ineps_nonlinear_figs:approx}). Below a given $\epsilon$-value, varying with the mesh size, the total error can be assimilated to $e_0$ and the decrease of $\epsilon$ is ineffective. The discrete operators being second order accurate $e_0$ is quadratically decreasing with the mesh step $h$.

As a consequence, we can conclude that $p_{\epsilon,app}$ converges to $p_{0}$ when $\epsilon$ converges to 0 alongwith $h$. This is exactly the Asymptotic-Preserving property of the scheme we intended to validate.

\section{Conclusions and perspectives}

In this paper we have presented an Asymptotic-Preserving numerical method for singular perturbation of non-linear anisotropic reaction-diffusion problems. The Asymptotic-Preserving property of the scheme is ensured thanks to a solution decomposition explained in full details in the most simple framework of a linear problem. This method is then generalized to non-linear problems thanks to Gummel's linearization method.

\indent In a second part, several two-dimensional numerical investigations of the AP-scheme are performed. These tests reveal a very weak dependence of the scheme accuracy with respect to the anisotropy direction, demonstrating the relevance of the use of non-adapted coordinates. The Asymptotic-Preserving property of the scheme is also validated for vanishing $\epsilon$ on linear as well as non-linear problems. The solution of the limit problem is accurately captured with no restrictions on the anisotropy strength. Furthermore, the computational efficiency of the method, in terms of memory as well as CPU usage, does not depend on this anisotropy strength. 

\indent Several applications of the present work can be investigated: at present time, the method has been used for the resolution of linear anisotropic diffusion problems for a two-fluid Euler-Lorentz model (see \cite{Brull-Degond-Deluzet-Mouton}) and the non-linear version of the method will be coupled to an Asymptotic-Preserving scheme for a one-fluid full Euler-Lorentz model (see \cite{Degond-Deluzet-Mouton}). %\textit{Other work pathes can be considered: indeed, we can consider an anisotropic diffusion problem in which there are some second order derivatives in the direction of $\mathbf{b}$ and in the orthogonal direction to $\mathbf{b}$ (see [\refcite{Narski}]).}

% \section{Introduction}\label{intro}
% 
% Put a general  introduction to your paper here. Separate text
% sections with other sections.
% 
% \section{Put title of the next section here}\label{an apprpriate
% label}
% 
%           %If you have subsections use:
% \subsection{Subsection title}\label{another label}
% 
% Don't forget to give each section, subsection, equation, theorem,
% corollary, etc. a unique label, and when you refer to the results
% later in the text use \ref{<labelname>} instead of explicitly
% writing the number of the environment in question.
% 
% This use of \ label and \ ref is REQUIRED for  papers.
% 
% Similarly, always use \cite{biblabelname} to refer to
% bibliographic references, which would then be entered in the
% bibliography via
%           %\bibitem{biblabelname}.
% 
%           %
%           % For figures use
% 
%           %\begin{figure}
% 
%           %The use of .eps files is encouraged, in which case you should
%           %un-comment the \uspackage{graphics} command above, and use the
%           %command
%           %\include{figure.eps}
%           % to insert the figure file.
% 
%           %\end{figure}
% 
% 
%           % BibTeX users please use
% 
%           % \bibliographystyle{}
% 
%           % \bibliography{}
% 
%           %
% 
%           % Non-BibTeX users please use
\medskip

{\bf Acknowledgement.} This work has been supported by the french magnetic fusion programme FR-FCM, by the INRIA large-scale initiative 'FUSION', by 'BOOST' and 'IODISSEE' ANR projects and by the CEA-Cadarache in the frame of the contract 'APPLA' (\# V3629.001 av. 2). The authors wish to thank P. Degond for suggesting this problem and for very fruitful discussions on the topic.

\end{document}